\input amstex
\input epsf
\documentstyle{amsppt}
\NoBlackBoxes
\define\R{\Bbb{R}}
\define\0h{\hat 0}
\define\1h{\hat 1}
\define\nd{\binom{[n]}d}
\define\In{\operatorname{In}}

\define\Max{\operatorname{Max}}
\define\Min{\operatorname{Min}}
\define\precdot{<\!\!\!\cdot\,\,}

\define\pr{\operatorname{Pr}}
\define\std{\operatorname{std}}
\define\sign{\operatorname{sign}}
\define\overbar{\overline}
\define\lk{\operatorname{lk}}
\define\adjoin{\amalg}

\topmatter 

\title Maps between higher Bruhat orders and higher Stasheff-Tamari posets 
\endtitle
\rightheadtext{Higher Bruhat orders}
\author Hugh Thomas \endauthor
\date December 3, 2002 \enddate

\abstract
We make explicit a description 
in terms of convex geometry of the higher Bruhat orders $B(n,d)$ 
sketched by Kapranov and Voevodsky.  
We give an analogous
description of the higher Stasheff-Tamari poset $S_1(n,d)$.  
We show that the map $f$ sketched  by Kapranov and Voevodsky from $B(n,d)$ 
to $S_1([0,n+1],d+1)$ coincides with the map 
constructed by Rambau, and 
is a surjection for $d \leq 2$.  We also give geometric descriptions of 
$\lk_{0}\circ f$ and $\lk_{\{0,n+1\}} \circ f$.  
We construct a map analogous to $f$ from
$S_1(n,d)$ to $B(n-1,d)$, and show that it is always a poset embedding.  
We also give an explicit criterion to determine if an element of $B(n-1,d)$
is in the image of this map.  

\endabstract
\endtopmatter

\document

\head 1. Introduction \endhead

The higher Bruhat orders $B(n,d)$ 
 were introduced by Manin and Schechtman [MS] in 
connection to discriminental hyperplane arrangements.  
They give a  combinatorial definition of $B(n,d)$ which
we shall review in the next section. The choice of name stems from the
fact that $B(n,1)$ is isomorphic to weak Bruhat order on the
symmetric group.  

Shortly following the definition of the higher Bruhat orders, Kapranov
and Voevodsky wrote a paper [KV] which presented two alternative 
interpretations for the higher Bruhat orders, in terms of oriented matroids,
and in terms of convex geometry.  The oriented matroid approach was 
later taken up by Ziegler [Zi].  The convex geometric approach has not been
significantly written about since.  It is the focus of the first part 
of our paper.  

The convex-geometric approach to the higher Bruhat orders is as follows.  
Consider the $n$-cube
$[-1,1]^n$.  Let $B(n,0)$ denote the set of 
vertices of the cube, with the usual Cartesian product order, so that 
$B(n,0)$ is a Boolean lattice.  Its minimum element is $(-1,\dots,-1)$,
and its maximum element is $(1,\dots,1)$.  Now let $B(n,1)$ be the set of 
increasing paths along edges of the cube 
from $(-1,\dots,-1)$ to $(1,\dots,1)$.  There are 
$n!$ of these, and they are naturally in bijection with $S_n$.  We put an
order on this set by defining covering relations: $\sigma \gtrdot \tau$ if
$\sigma$ and $\tau$ coincide except on the boundary of some 2-face, where
$\sigma$ uses the top two edges of the face, and $\tau$ uses the bottom
two edges.  (We will be more explicit about how to understand ``bottom'' and
``top'' in the next section.)  Under this order, the poset $B(n,1)$ is
isomorphic to weak Bruhat order on the symmetric group.  

Now consider  
collections of 2-faces of the cube which form a homotopy from the minimum path
to the maximum path, 
and which
are non-backtracking.  (This non-backtracking condition generalizes the
``increasing'' condition in the dimension 1 case. We shall give more 
precise definitions 
in the next section.)  These homotopies form the 
elements of $B(n,2)$.  As before, the order on $B(n,2)$ is defined by
specifying covering relations: $\sigma \gtrdot \tau$ if $\sigma$ and 
$\tau$ coincide except on the boundary of a 3-face, where $\sigma$ 
uses the top three facets, and $\tau$ uses the bottom three facets.  The other
$B(n,d)$ are defined similarly.  
The first goal of this paper is to write down this description explicitly,
and to show that it is equivalent to the combinatorial definition of
[MS]. 

In order to describe the second goal of the paper, we must now turn to the
higher Stasheff-Tamari posets.  In fact, there are two different 
Stasheff-Tamari posets structures $S_1(n,d)$ and $S_2(n,d)$ defined on the
same set of objects $S(n,d)$.  We shall only be interested in the first of
these posets, so we shall suppress the subscript.  $S(n,d)$ is usually 
viewed as the set of triangulations of the cyclic polytope $C(n,d)$. 
We will give an equivalent definition, 
analogous to the one above for higher Bruhat orders where the cube has
been replaced by a simplex.  

In order to define $S(n,d)$,
start with an $n-1$-simplex, with 
vertices labelled from $1$ to $n$.  $S(n,0)$ is the set of vertices,
with the order given by the labelling.  
The objects of $S(n,1)$ are the increasing paths from the bottom
vertex to the top vertex.  The order is by reverse refinement: 
the bottom path is the path that includes 
every vertex, and the top path is the one that uses only $1$ and $n$.  
The objects of $S(n,2)$ are the sets of 2-faces of the simplex which form a 
non-backtracking homotopy from the bottom path to the top path.  The order
on $S(n,2)$ is defined by specifying covering relations: $S \gtrdot T$ if
$S$ and $T$ coincide except on the boundary of a 3-simplex, where $S$ 
contains the upper faces and $T$ the lower faces.  The higher $S(n,d)$ are
defined similarly.

In [KV], a map called $f$ from $B(n,d)$ to 
$S([0,n+1],d+1)$ was described as follows.  (We write $S([0,n+1],d+1)$ 
to indicate that the vertices are labelled by the numbers 
from 0 to $n+1$.)  
There is a poset isomorphism from vertices of
the $n$-cube (i.e. $B(n,0)$) to elements of $S([0,n+1],1)$: namely, the 
coordinates of the
vertex which are negative tell you which vertices belong in the path
in addition to $0$ and $n+1$.  Now, an element of $B(n,1)$, which
is a path through the $n$-cube, determines a sequence of vertices of the cube. 
We apply the
map from $B(n,0)$ to $S([0,n+1],1)$ 
to each vertex in succession, to get a sequence of paths through the 
$n+1$-simplex.  Two successive paths differ in that one vertex which is present
in the first path is not present in the second.  To each pair of successive
paths, we associate the triangle whose vertices are the removed vertex and
its two neighbours along the path.  
These triangles form a homotopy from the bottom path 
(which contains all the vertices)
to the top path (which contains only the end-points), and hence define an
element of $S([0,n+1],2)$.  In the example below, the bold path through the
cube on the left gives rise to the triangulation shown on the right.  

$$\epsfbox{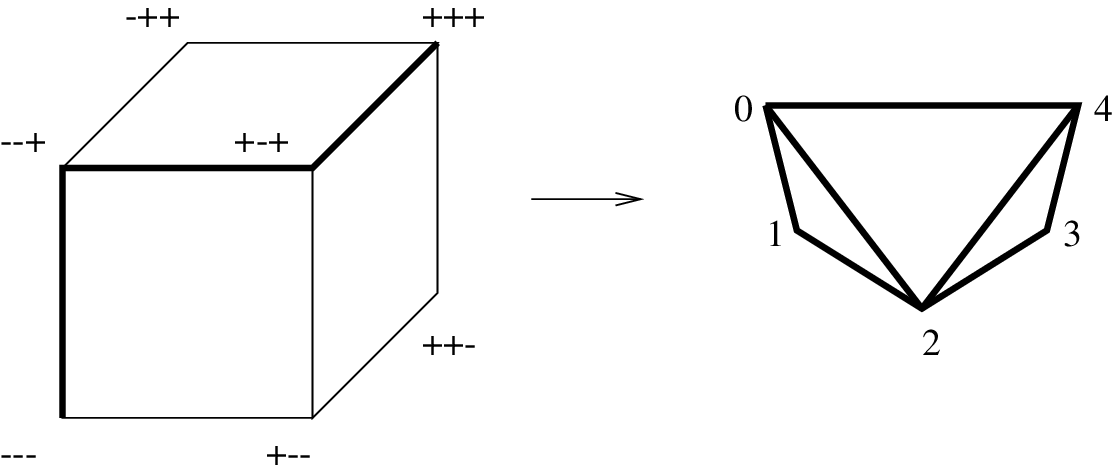}$$

It is claimed in [KV]
that one can define a similar map 
$f:B(n,d)\rightarrow S([0,n+1],d+1)$
for all $d$, and further that this map
is surjective for all $d$.  Rambau [Ra1] constructed an explicit map from 
$B(n,d)$ to $S([0,n+1],d+1)$, but he did not show that it coincided with
the map described in [KV].  In the second part of the paper, we prove
that the map described in [KV] does indeed coincide with that defined in [Ra1],
and that it is surjective for $d \leq 2$.  

In the second part of the paper, we also give
geometric interpretations of two maps associated to $f$.  
For $S\in S([0,n+1],d+1)$, $\lk_0 (S)$ is the link of $S$ at zero,
which can be viewed in a natural way as lying in $S(n+1,d)$.  
We can also take $\lk_{\{0,n+1\}}(S)$; this lies in $S(n,d-1)$.  
We show that, for $\pi \in B(n,d)$, $\lk_{\{0,n+1\}}(f(\pi))$ 
coincides with the vertex figure of $\pi$ at 
$(1,\dots,1) \in [-1,1]^n$. 
We also give a similar interpretation for 
$\lk_0(f(\pi))$.

The third part of our paper consists of the construction of a map 
$g:S(n,d) \rightarrow B(n-1,d)$.  
As already explained, $S(n,0)$ is a chain of $n$ elements which we view
as the vertices of an $n-1$-simplex.  We map vertex $a$ to the corner
of $[-1,1]^{n-1}$ whose final $a-1$ coordinates are +1, and the others -1.  
An element of $S(n,1)$ is mapped to an increasing path through the
$n-1$-cube which passes through the vertices corresponding to the 
vertices on the path through the $n-1$-simplex.  For $S \in S(n,2)$, 
$g(S)$ is defined as the unique homotopy from the minimal path through
the cube 
to the 
maximal path through the cube which passes through all the paths corresponding to paths 
through the $n-1$-simplex along edges in $S$.  An analogous statement
holds for $d>2$.  

We give several equivalent definitions for $g$, including two which are
 explicit and non-inductive.
  We also show that the  map $g$ is a 
poset embedding, and we give an explicit criterion to determine if an
element of $B(n-1,d)$ is in the image of $g$.  This amounts to giving a 
new equivalent definition of $S(n,d)$ without reference to convex geometry.

\head 2. Higher Bruhat orders \endhead

We begin by recalling the definition of the higher 
Bruhat order $B(n,d)$ for $d 
\leq n$
 positive integers, given by Manin and Schechtman [MS].  
We write $\binom{[n]}d$ for the set of subsets of $[n]$ of size $d$.  
A $d$-packet consists of the subsets of size $d$ of a set of $d+1$ integers.
A total order on $\binom{[n]}d$ is admissible if every $d$-packet occurs in
either lexicographic order or its reverse. The set of admissible
orders on $\binom{[n]}d$ is called $A(n,d)$.  

Two admissible orders are said to be equivalent
if they differ by a sequence of transpositions of adjacent elements not
both lying simultaneously in any $d$-packet.  If $\pi$ is an admissible order,
we write $[\pi]$ for its equivalence class.  

$B(n,d)$ is a poset whose elements are the equivalence classes of admissible
orders on $\binom{[n]}d$.  The order is given by specifying covering
relations. Let $\pi \in A(n,d)$.  Suppose there is some $d$-packet which 
occurs consecutively in $\pi$, and in lexicographic order.  Let $\sigma$
denote the (automatically 
admissible) order obtained by reversing this $d$-packet.  Then
$[\pi]\lessdot[\sigma]$ in $B(n,d)$.  The order on $B(n,d)$ is the transitive
closure of these covering relations.  

There are two orders on $\nd$ which are clearly admissible.  
Let $\0h_d$, $\1h_d$ denote the class in $B(n,d)$ of the lexicographic 
order and its reverse respectively.  It is clear that these elements are
minimal and maximal respectively in $B(n,d)$; in fact, they are its minimum and
maximum elements, see [MS].

There is a map $I: A(n,d) \rightarrow \Cal{P}(\binom{[n]}{d+1})$ which 
associates to any $\pi \in A(n,d)$ the set of $d$-packets which occur in
reverse order.  $I(\pi)$ is called the inversion set of $\pi$.  This
generalizes the usual notion of inversion set for a permutation.  
It is clear that $I$ is constant on equivalence classes, and so passes to 
$B(n,d)$.  As a map from $B(n,d)$, $I$ is injective.  

A subset of $\binom{[n]}{d+1}$ is said to be consistent if its restriction
to any $(d+1)$-packet consists of either an initial or a final subset
with respect to lex order.  
Ziegler showed in [Zi] that a subset
of $\binom{[n]}{d+1}$ is in the image of $I$ iff it is consistent.

It will be convenient to define $B(n,0)$ to be the set of subsets of $[n]$,
ordered by inclusion.  
The inversion set of an element of $B(n,0)$ is just
the set itself.  $\1h_0 = [n]$; $\0h_0=\emptyset$.  

We now give the convex-geometric definition of 
the higher Bruhat orders, formalizing
ideas from [KV].   
%To each element 
%$[\pi] \in B(n,d)$ we will associate a set $K(\pi)$ of 
%$d$-faces of an $n$-cube.  The 
%union of these $d$-faces will be homeomorphic to a $d$-dimensional disk,
%whose boundary does not depend on $\pi$.  For $d=1$, this will 
%restrict to the well-known correspondence between $S_n$ and 
%paths of length $n$ from a fixed vertex of an $n$-cube to its opposite
%vertex.
$[-1,1]^n$  will be our standard $n$-cube.  We shall
keep track of its faces as maps from $[n]$ to the
set $\{-1,*,1\}$, where a $d$-face will have $*$ occurring in $d$ places,
these being the dimensions in which the face extends.  We will sometimes
refer to a set of faces of the cube when what we mean is the union of the
set of faces.

For $G$ a set, let $\Xi_G(x)=1$ if $x\in G$ and $-1$ otherwise.  
For $X=\{a_1,\dots,a_d\}_< \subset [n]$, $y \not \in X$, define:
$$p(y,X)=\left\{ \matrix 1 & y<a_1 \\ (-1)^i & a_i < y < a_{i+1} \\
(-1)^d & a_d<y \endmatrix \right. $$

Fix $\alpha \in B(n,d)$.  
For each $X \in \nd$, let
$$F^\alpha_X(i)=\left\{\matrix * & i \in X \\
                          p(i,X)\Xi_{I(\alpha)}(X\cup
\{i\}) & i \not \in X
\endmatrix \right. $$

Let $K(\alpha)$ consist of the $F^\alpha_X$ for all $X$.  
%For convenience, let $c^\alpha(i,X)=p(i,X)\Xi_{I(\alpha)}(X\cup
%\{i\})$.

We identify linear maps from $\R^n$ to $\R^d$ with $d \times n$ matrices.  
We say that a map $T: \R^n \rightarrow \R^d$ is totally positive if 
the determinants of all its
minors are positive.  (Note that there are many totally positive matrices, for
example, a Vandermonde matrix $T_{ij}=c_j^i$ with $0<c_1<\dots<c_n$.  The
determinant of any minor of this matrix equals a Vandermonde determinant times
a Schur function, both of which are positive.)

We say that a collection of convex sets tiles a region if the region is the
union of the convex sets and 
the sets overlap only on boundaries.  

The main theorem of this section is the following:

\proclaim{Theorem 2.1} For any $\alpha \in B(n,d)$, the set $K(\alpha)$ of
$d$-faces of the standard $n$-cube is homeomorphic to a disk, has boundary
$K(\0h_{d-1})\cup K(\1h_{d-1})$, and the image of the $K(\alpha)$ under any 
totally positive map $T$ from $\R^n$ to $\R^d$ forms a tiling of the 
image of the standard $n$-cube under $T$.

Conversely, given a set $K$ of $d$-faces of the standard $n$-cube, 
such that the images under some totally positive map $T$ 
of the faces in $K$ tile the image of the standard $n$-cube, 
it follows that $K=K(\alpha)$ for some $\alpha \in B(n,d)$.  \endproclaim

\demo{Proof}
Given a convex polytope in $\R^d$, we say that a facet is an upper facet if
the polytope lies below it with respect to the final co-ordinate, and similarly
for lower facets.  A facet parallel to $e_d$ is neither upper nor lower, but
the cases in which we shall be interested will exclude that possibility,
so that every facet is either upper or lower.  We now prove a simple lemma
about upper and lower facets of images of cubes:

\proclaim{Lemma 2.1} Let $W$ be a totally positive map from $\R^{d}$ to 
$\R^d$.  The upper facets of the image under $W$ of $[-1,1]^d$ are 
the $T_i$ defined by 
$$T_i(j)=\left\{ \matrix * & j \ne i \\
                       (-1)^{d+j} & j=i \endmatrix \right. , $$
while the lower facets are the $L_i$ defined by   
$$L_i(j)=\left\{ \matrix * & j \ne i \\
                       (-1)^{d+j+1} & j=i \endmatrix \right. . $$\endproclaim

\demo{Proof} Since $W$ is totally positive, its inverse satisfies 
$\sign((W^{-1})_{ij})=(-1)^{i+j}$.  Thus,
the inverse image under 
$W$ of $(0,0,\dots,0,1)$ will be alternating in sign, with its last entry
positive, and the desired result follows. \enddemo

Observe that if $T$ is a totally positive map from $\R^n$ to $\R^d$, then
the restriction of $T$ to any $d$-dimensional co-ordinate subspace of 
$\R^n$ is a map to which Lemma 2.1 applies.  Thus, Lemma 2.1 tells us about the
upper and lower facets of the image under $T$ of any $d$-dimensional face
of $[-1,1]^n$.   We shall sometimes speak of the upper or lower facets of
a face of $[-1,1]^n$ when what we mean is facets whose images are upper or 
lower in the image of the face under any totally positive map.

We now begin to prove the forward direction of the theorem.  
The proof is by induction on $d$.  It is clear for $d=0$. Assume
it holds for all dimensions less than $d$.

Let $\pi \in A(n,d)$.  Let $Y_1,\dots,Y_{\binom{n}{d}}$ be the elements of
$\binom{[n]}d$ under the order $\pi$.  Let
$\In_i(\pi)$ denote $\{Y_1,\dots,Y_i\}$.  
Observe that for $0 \leq i \leq \binom
nd$, $\In_i(\pi)$ is consistent.  
Thus, we may define a sequence $\alpha_j \in B(n,d-1)$ by
$I(\alpha_j)=\In_i(\pi)$.  
Note that $\alpha_0=\0h_{d-1}$; $\alpha_\binom nd
=\1h_{d-1}$.  

Let $\pr$ denote the map from $\R^d$ to $\R^{d-1}$ forgetting the last 
coordinate.  Let $T'=\pr\circ T$.  Then $T'$ is totally positive, so 
by the induction hypothesis, 
for $0 \leq j \leq \binom nd$, the images under $T'$ of 
$K(\alpha_j)$ define tilings of $T'([-1,1]^n)$.  Let $\Gamma_j$ be the
image under $T$ of $K(\alpha_j)$. Then for each $x \in T'([-1,1]^n)$,
$\Gamma_j \cap \pr^{-1}(x)$ consists of a single point.

Fix $0 \leq i < \binom nd$. 
For $Z \in \binom{[n]}{d-1}$,  
$F_Z^{\alpha_i}$ and $F_Z^{\alpha_{i+1}}$ are
the corresponding faces of $K(\alpha_i)$ and $K(\alpha_{i+1})$.  
$F_Z^{\alpha_i}$ and $F_Z^{\alpha_{i+1}}$ coincide except for 
$Z \subset Y_{i+1}$.  
The $d$ faces of each not shared by the other are the $2d$ faces of 
$F_{Y_{i+1}}^{[\pi]}$.  So $\Gamma_i$ and $\Gamma_{i+1}$ coincide except that 
each of them contains $d$ of the $2d$ faces of $T(F_{Y_{i+1}}^{[\pi]})$.  
One checks using Lemma 2.1 that 
$\Gamma_i$ 
contains the lower faces and $\Gamma_{i+1}$ contains the upper faces.  

Thus it follows that if $i<j$, for any point $x\in T'([-1,1]^n)$, the 
intersection of $\Gamma_i$ with $\pr^{-1}(x)$ lies on or 
below the intersection of
$\Gamma_j$ with $\pr^{-1}(x)$.  

Thus, the images of the $K([\pi])$ intersect only on boundaries, since they 
are 
separated by the $\Gamma_i$, and they therefore tile the region between
the images of $K(\0h_{d-1})$ and $K(\1h_{d-1})$.  
Also, since this holds for any $\pi$, it follows
that the region between the images of $K(\0h_{d-1})$ and $K(\1h_{d-1})$ is 
the entire image of
$[-1,1]^n$ under $T$, as desired.  It also follows from this that the images
of the $K(\0h_{d-1})$ are exactly the bottom facets of the image under $T$ 
of $[-1,1]^n$, while the images of the $K(\1h_{d-1})$ are its top facets. 

Another result of what we have shown so far is that every face in every 
$K(\alpha_i)$ occurs as a facet of either one or two faces in $K([\pi])$:
one if the face we are interested in is $K(\0h_{d-1})$ or $K(\1h_{d-1})$ and
two otherwise.  This allows us to conclude that $T$ restricted to 
$K([\pi])$ is a homeomorphism to the image of $[-1,1]^n$, which is clearly
(homeomorphic to) a disk.

We now turn to the converse direction of the theorem.  The proof is again
by induction on $d$.  Again, it is obvious for $d=0$, so we assume
$d>0$, and that the converse holds for dimension $d-1$. 

Fix $K$ a set of
$d$-faces of the standard $n$-cube, 
as in the statement of the theorem.  We will
now define a sequence $\alpha_0,\dots,\alpha_{\binom nd}$ of elements of 
$B(n,d-1)$ such that any face in $K(\alpha_i)$ for any $i$ occurs as a facet
of some face in $K$.  

Let $\alpha_0=\0h_{d-1}$.  Inductively, given $\alpha_i$, for $i <\binom nd$,
we will define
$\alpha_{i+1}$ as follows.  Let $\Gamma_i$ be the image of $K(\alpha_i)$ 
under $T$.  $\Gamma_i$ divides the faces of $K$ into those whose images 
are above $\Gamma_i$ and those below it.  I claim that there exists some 
$d$-face in $K$ all of whose lower facets are in $K(\alpha_i)$.
We will pick one such, and call it $Y_{i+1}$.  To find such a $Y_{i+1}$,
pick any $X \in K$ which lies above $\Gamma_i$.  If $X$ has some lower facet
not in $K(\alpha_i)$, replace $X$ by the face in $K$ which contains
this lower facet of $X$ as an upper facet.  This new face still has its 
image under $T$ lying above $\Gamma_i$, but its topmost point is lower than
that of the old face.  Thus, this process cannot loop back on itself, but 
must terminate, and it must terminate in a face $Y_{i+1}$ 
which has all its lower
facets in $K(\alpha_i)$.  
(A similar statement for cyclic polytopes is proved in [Ra1].)  

Now, let $J$ denote $K(\alpha_i)$ with the lower facets of $Y$ replaced
by its upper facets.  It is clear that $J$ satisfies the conditions of 
the theorem, and thus that, by induction, 
we can define $\alpha_{i+1}$ by saying that
$K(\alpha_{i+1})=J$.  We check that 
$I(\alpha_{i+1})=I(\alpha_i)\cup\{Y_{i+1}\}$.

Now let $\pi$ denote the order on $\binom{[n]}d$ given by 
$Y_1,\dots,Y_{\binom nd}$. Since $\In_i(\pi)=I(\alpha_i)$, and 
the $I(\alpha_i)$ are all consistent, it follows that 
$\pi$ is an admissible order.  

Finally, we observe that $K([\pi])=K$ because both $K$ and $K([\pi])$ can
be characterized as the set of faces between $\alpha_i$ and $\alpha_{i+1}$
for some $i$.  This completes the proof of Theorem 2.1.  \enddemo

We have shown that the elements of $B(n,d)$ can be represented as sets of
$d$-faces of the standard $n$-cube.  To describe the order relation 
on $B(n,d)$ in terms of this description, we have the following proposition:

\proclaim{Proposition 2.1} For $\sigma, \tau \in B(n,d)$, $\sigma \gtrdot \tau$
iff $K(\sigma)$ and $K(\tau)$ coincide except on the facets of a 
$d+1$-cube, where $K(\sigma)$ contains the upper facets and $K(\tau)$ 
contains the lower facets.  \endproclaim

\demo{Proof} Suppose that $\sigma \gtrdot \tau$.  So $I(\sigma)=I(\tau)
\cup \{Y\}$. 
Choose an admissible order $\pi$ on $\binom{[n]}{d+1}$ so that
$I(\tau)$ preceds $Y$ precedes the rest of $\binom{[n]}{d+1}$.  Now,
as in the proof of Theorem 2.1, we see that $\sigma$ and $\tau$ differ only
on the boundary of $F^{[\pi]}_Y$: $\sigma$ contains its upper facets,
and $\tau$ its lower facets, as desired.  

To prove the converse, observe that since $K(\sigma)$ and $K(\tau)$ 
coincide outside the boundary of a $d+1$-face, say $Y$, $I(\sigma)$ and
$I(\tau)$ must coincide except as regards containment of $Y$, and the
desired result follows.  \enddemo

\head 3. The higher Stasheff-Tamari posets  \endhead

As explained in the introduction, the usual way of thinking of the higher
Stasheff-Tamari posets $S(n,d)$ is as a poset on the set of triangulations of
a cyclic polytope.  
To motivate the existence of a connection
to the higher Bruhat orders, a different definition, one analogous to the
convex-geometric definition of $B(n,d)$ given above, will be more relevant.
To avoid confusion, we shall give the poset we define in this manner a new
name, $T(n,d)$, and then prove that $T(n,d)$ coincides with $S(n,d)$.  

The standard $n-1$-simplex, $\Delta_{n-1}$, is the convex hull of the basis
vectors in $\R^n$.  Its $d$-faces are indexed by $d+1$-subsets of $[n]$,
which designate which vertices lie on the face.

If $W:\R^n \rightarrow \R^d$, let $\overbar W$ be the linear map from 
$\R^n$ to 
$\R^{d+1}$ defined by setting $\overbar W(e_i)=(1,W(e_i))$.  In terms of
matrices, we can say that the matrix for $\overbar W$ is obtained from
that for $W$ by adding a first row of all ones.  
We say that 
$W$ is {\it affinely positive} if $\overbar W$ is totally positive. 

We now prove a lemma analogous to Lemma 2.1, but concerning simplices, not
cubes.  

\proclaim{Lemma 3.1}
Let $W$ be an affinely positive map from $\R^{d+1}$ to $\R^d$.  
The top facets of the image of $\Delta_d$ are those which omit a vertex
with the same parity as $d$;
the bottom facets are those which omit a vertex of opposite parity to 
$d$.  
\endproclaim

\demo{Proof} Consider the totally positive map $\overbar W:\R^{d+1}\rightarrow
\R^{d+1}$.  Let $x={\overbar W}^{-1}(e_{d+1})$.  Then $x$ is parallel
to the affine span of $\Delta_{n-1}$, and, as in the proof of Lemma 2.1, we 
are interested in whether $x$ points into or out of each facet of 
$\Delta_{n-1}$.  This is equivalent to asking whether $x$ points into or out
of the corresponding facet of the cone over $\Delta_{n-1}$ with cone point the
origin, and now we use the fact that, as in Lemma 2.1, $x$ alternates in sign. 
\enddemo

Since this does not depend on the choice of $W$, this allows us to refer
to ``upper'' or  ``lower'' facets of a face of $\Delta_{n-1}$, meaning
facets whose images are upper or lower in the image of the face under any 
affinely positive map.  

We can now define $T(n,d)$.  

\proclaim{Definition of $T(n,d)$}
An element $S$ of $T(n,d)$ is a set of $d$-faces of $\Delta_{n-1}$, with the
property that under some  affinely positive map $W$, 
the images under $W$ of the faces in $S$ tile                        
$W(\Delta_{n-1})$.  The order on $T(n,d)$ is defined by giving covering 
relations:
$S\gtrdot T$ iff $S$ and $T$ coincide except on the boundary of a 
$d+1$-simplex, 
where $S$ contains the upper facets of the simplex and $T$ contains the lower
facets.   \endproclaim

\proclaim{Proposition 3.1}
If $S \in T(n,d)$ then for any affinely positive map $V$, the images 
under $V$ of 
the faces in $S$ form a tiling of  $V(\Delta_{n-1})$.  \endproclaim

\demo{Proof}
First, we show how much information we need about a set of
points to be able to determine when a collection of simplices forms a 
tiling.  

\proclaim{Lemma 3.2} Given $n$ points in $\R^d$ no $d+1$ of which lie on any
affine hyperplane:

(i) The boundary facets of the convex hull of the set of vertices
are the $d$-sets of vertices having the property that all the other
vertices lie on the same side of their affine span.  

(ii) A collection of $d+1$-sets of vertices (``simplices'') 
forms a tiling iff every
facet of every simplex is either a boundary facet of the convex hull
and appears as a facet of exactly one simplex, or else appears
as a facet of exactly two simplices, and the vertices of these
two simplices not on the shared facet lie on opposite sides of the
affine span of the shared facet.  \endproclaim

\demo{Proof} Part (i) is obvious. 
To establish part (ii), let $P$ be the 
convex hull of the points.  Let $S$ be a set of simplices.  If $S$ is 
a tiling, it is clear that it has the above properties.  Now, 
we assume it has the above properties, and we wish to show that $S$ forms
a tiling.  

Pick an arbitrary direction vector $v$.  Now, pick a line with direction vector
$v$ which passes through $P$, and does not intersect any simplex in $S$
in a face of codimension more than $1$.  We say that a point on the line is
bad if it neither lies on the boundary of a simplex in $S$, 
 nor lies in exactly one simplex of $S$.  
Let $x$ be the point furthest
along the line in the closure of the set of bad points.

First, consider the case where $x$ is in the interior of $P$.
Since the points just beyond
$x$ are good, they lie in exactly one simplex, say $A$, of $S$, 
and clearly $x$ lies
on the boundary of $A$.  The facet of $A$ containing $x$ lies in
exactly one other simplex of $S$, say $B$.  Since the vertex of $B$ not 
lying on the facet containing $x$ lies on the opposite side from $A$, the
points immediately before $x$ lie in $B$.  This shows that the vertices
just before $x$ lie in at least one simplex of $S$.  Suppose they lie in
another one as well, say $C$.  Then $x$ must lie on the boundary of $C$, and,
as before, the points just past $x$ must lie in another simplex, say $D$.
But since the points just past $x$ were assumed to be good, this is impossible.
The case where $x$ is on the boundary of $P$ is similar.  

We have now showed that the points not on any boundary face of a simplex
of $S$, lying on a line in direction $v$ which
doesn't intersect any faces of $S$ in codimension more than 1, all lie in
exactly one simplex of $S$.  But this set is dense in $P$, so the simplices
of $S$ form a tiling, as desired.  This completes the 
proof of Lemma 3.2. \enddemo

Proposition 3.1 will now follow from Lemma 3.2 and the following lemma:

\proclaim{Lemma 3.3} Let $V$ be an affinely positive map from $\R^{n}$ to
$\R^{d}$.  Let $x_i=V(e_i)$.  Then no $d+1$ of the $x_i$ lie on a common
affine hyperplane, and for any $a_1,\dots,a_{d+1},i,j$ distinct integers 
in $[n]$, whether or not $x_i$ and $x_j$ lie on opposite sides of the
affine hyperplane spanned by $x_{a_1},\dots,x_{a_{d+1}}$ does not depend on
$V$.  \endproclaim

\demo{Proof}We begin with a lemma:

\proclaim{Lemma 3.4} Let $T$ be a totally positive map from $\R^{n}$ to 
$\R^{d+1}$.
Let $x$ be a non-zero vector in $\ker(T)$.  
Then $x$ has at least $d+2$ non-zero components.  For any set of $d+2$ 
components, there is an $x \in \ker(T)$ with exactly those components 
non-zero.   
If $x\in \ker(T)$ has exactly $d+2$ non-zero components, then its non-zero 
components
alternate in sign.  \endproclaim

\demo{Proof} Since $T$ is totally positive, its restriction to any 
co-ordinate subspace of dimension $d+1$ is non-singular, so no non-zero
element of any of these subspaces could be in the kernel of $T$.  Thus
$x \in \ker(T)$ implies that $x$ has at least $d+2$ non-zero components.  

The restriction of $T$ to any co-ordinate subspace of dimension $d+2$ 
must have a non-zero kernel, but if $x$ is a non-zero element of the kernel,
by what we have already shown, it must have all $d+2$ components non-zero.  
This shows that for any choice of $d+2$ components, there is an element of 
the kernel of $T$ with exactly those componenets non-zero.  

Assume $x$ has exactly $d+2$ non-zero components: 
$x=\sum_{i=1}^{d+2} c_ie_{a_i}$
with $a_1<\dots<a_{d+2}$.  Let $S$ be the restriction of $T$ to the span
of the $e_{a_i}$ for $1\leq i \leq d+1$.  Let $y=c_{d+2}T(e_{a_{d+2}})$.
Then $S^{-1}(y)=\sum_{i=1}^{d+1}c_ie_{a_i}$.  Computing $S^{-1}(y)$, we see
that its coefficients in the $e_{a_i}$ are determinants of minors of $T$, up 
to an alternating
sign, which proves the final statement.  \enddemo

We now return to the proof of Lemma 3.3.  
Let $V$ be an affinely positive map from $\R^n$ to $\R^d$.  Let
$x_k=V(e_k)$.  No $d+1$ of the $x_k$ lie on any affine hyperplane, since
this would imply a linear dependence among the corresponding $\overbar V(e_k)$,
which is impossible because $\overbar V$ is totally positive.  

Now suppose $x_i$ and $x_j$ lie on opposite sides of the
affine span of $x_{a_1},\dots,x_{a_{d}}$. Then there exists $0<c<1$ such that
$cx_i +(1-c)x_j$ lies in the affine span of $x_{a_1},\dots,x_{a_{d}}$, or in 
other words that there is some $0<c<1$ and some $b_1,\dots,b_d$ summing to 1
such that $cx_i + (1-c)x_j = \sum b_kx_{a_k}$.  Thus, 
$ce_i +(1-c)e_j - \sum b_ke_{a_k}$ lies in the kernel, not merely of $W$,
but in fact of $\overbar W$. By Lemma 3.4, we know that this implies that
there are an even number of $a_k$ lying between $i$ and $j$.  Conversely,
if there are an even number of $a_k$ lying between $i$ and $j$, we can
reverse the argument to show that $x_i$ and $x_j$ lie on opposite sides
of the affine span of the $x_{a_k}$.  Thus, we see that whether or not
$x_i$ and $x_j$ lie on opposite sides of the affine span of 
$x_{a_1},\dots,x_{a_{d}}$ does not depend on the choice of $V$.
This completes the proof of Lemma 3.3 (and hence also of Proposition 3.1). 
\enddemo \enddemo

We can now define two elements of $T(n,d)$, $\1h_d$ and $\0h_d$, as follows.
Let $W$ be an affinely positive map from $\R^n$ to $\R^d$.  As we have
already shown (Lemma 3.2 (i)), the boundary
facets of the the image of $\Delta_n$ under $W$ do not depend on $W$, so 
let
$\1h_d$ consist of the faces of $\Delta_{n-1}$ corresponding to 
upper boundary facets of $W(\Delta_{n-1})$, and let $\0h_d$ consist of
its lower boundary facets. We remark that $\1h_d$ is clearly a maximal 
element of $T(n,d)$, and $\0h_d$ is clearly a minimal element.  They are
in fact maximum and minimum respectively, which we know from [Ra1]   
(once we know that
$T(n,d)$ coincides with $S(n,d)$).  

The following proposition is now clear:

\proclaim{Proposition 3.2} If $S \in T(n,d)$, 
the faces in $S$ are homeomorphic
to a disk, and their boundary equals $\0h_{d-1} \cup \1h_{d-1}$.  
\endproclaim

Finally, we show that $T(n,d)$ coincides with the poset $S(n,d)$ as 
conventionally defined.  We begin by reviewing the definition of $S(n,d)$.

Fix $d$ a positive integer. Let  $M(t)=(t,t^2,\dots,t^d)$.  Choose $n$
real numbers $t_1<\dots<t_n$.  Let $P$ be the convex hull of the $M(t_i)$.
$P$ is called a cyclic polytope.  

Many combinatorial properties of $P$ depend only on $d$ and $n$, and not
on the choice of $t_i$. Let $I$ be a $d$-set contained in $[n]$. Then
whether or not the $M(t_i)$ for $i \in I$  form a boundary facet of $P$
does not depend on the choice of $i$.  (In fact, the boundary facets 
are described
by the well-known ``Gale's Evenness Criterion,'' see [Gr].)  Further,
let $S \subset \binom {[n]}{d+1}$.  
To each $A\in S$, we can associate a simplex contained in $P$.  And again,
whether or not this collection of simplices forms a triangulation of $P$ 
does not depend on the choice of the $t_i$.  
Thus, we shall usually refer to ``the'' cyclic polytope in dimension $d$
with $n$ vertices, and denote it $C(n,d)$.  When we wish to emphasize
the choice of some particular $t_i$, we speak of a geometric realization of
$C(n,d)$. 

The partial order on $S(n,d)$ is given by describing its covering relations.
If one is familiar with the language of bistellar flips, one can say 
that the covering relations $S \precdot T$ are given by pairs $S$ and $T$ which
are related by a single bistellar flip, where bistellar flips are given
a certain natural orientation to determine whether $S$ precedes $T$ or
vice versa.  
The reader interested in a thorough explanation of this can consult [ER].

More explicitly, we can define the covering relations as follows, following
[Ra1].  
Let $M'(t)=(t,t^2,\dots,t^{d+1})\in \Bbb{R}^{d+1}$.  Pick
$t_1<\dots<t_n$.  This yields geometric realizations of $C(n,d+1)$ and
$C(n,d)$, where the map forgetting the last co-ordinate maps $C(n,d+1)$ down
to $C(n,d)$.  A triangulation $S\in S(n,d)$ defines a section $\Gamma_S \subset
\Bbb{R}^{d+1}$ over $C(n,d)$
by lifting its vertices $M(t_i)$ to $M'(t_i)$ and then extending linearly
over the simplices of $S$.   
Now, $S \precdot T$ precisely if $\Gamma_S$ and $\Gamma_T$ coincide 
except within the convex hull of $d+2$ vertices, where $\Gamma_S$ forms
the bottom facets and $\Gamma_T$ the top facets of a $d+1$-dimensional
simplex.  

Now we are ready to prove the following proposition:

\proclaim{Proposition 3.3} The poset $T(n,d)$ and the poset $S(n,d)$ coincide.  
\endproclaim

\demo{Proof}
Choose some $0<t_1<\dots<t_n$.  
If we define
a map $W$ from $\R^n$ to $\R^d$ by setting
$W_{ij}=t_j^i$, then $W$ is affinely positive (since
the determinants of its minors are given by a Schur function times a 
Vandermonde determinant, both of which are positive), and the image under
$W$ of $\Delta_{n-1}$ is exactly the geometric realization of $C(n,d)$ with
parameters $t_1,\dots,t_n$.  
Thus, by what we have already proven, the elements of $T(n,d)$ are 
in one-to-one correspondence with tilings of $C(n,d)$ by simplices whose
vertices are among the vertices of $C(n,d)$.  In general, a tiling by
simplices is not necessarily a triangulation, but because no $d+1$ of the
vertices of $C(n,d)$ lie on an affine hyperplane, the two notions 
coincide.  

It is also
easy to see that the covering relations in the two partially ordered sets
coincide.  This completes the proof of Proposition 3.3.
\enddemo

Since we have shown that $T(n,d)$ and $S(n,d)$ coincide, we shall use the
conventional notation of $S(n,d)$, but the reader is advised that we will
tacitly use the intuition that the elements of $S(n,d)$ can be considered
as sets of $d$-faces of an $n-1$-simplex.  

\head 4. The map from $B(n,d)$ to $S([0,n+1],d+1)$ \endhead

We will now proceed to elucidate the poset map sketched in [KV] from 
$B(n,d)$ to $S([0,n+1],d+1)$.  The definition is by induction. 

\proclaim{Definition 1} 
If $\alpha \in B(n,0)$, then $f_1(\alpha)$ is the path whose vertices 
are the elements of $[0,n+1]$ not in
$\alpha$, in increasing order.  

For $d>0$, let $\pi \in A(n,d)$.  Define $\alpha_i$ by 
$I(\alpha_i)=\In_i(\pi)$.  Let $S_i=f(\alpha_i) \in S([0,n+1],d)$.  We 
claim that for all $i$, either $S_i$ and $S_{i+1}$ coincide, or they differ
precisely
in that $S_i$ contains the bottom facets of some $d+1$-simplex $A_{i+1}$, 
while $S_{i+1}$
contains its top facets.  
Then, $f_1([\pi])$ consists of 
the collection of the $A_{i+1}$ for all $i$ for which $S_i$ and $S_{i+1}$
are different.  \endproclaim

Because of the reliance on the claim mentioned, this definition doesn't 
establish the existence of $f_1$.  We shall now give an explicit 
definition of 
$f_2$, essentially the map called $\Cal{T}_{\text{dir}}$ in [Ra1].  An induction
argument will then show that $f_2$ satisfies Definition 1.  

If $X=\{a_1,\dots,a_{d}\}_< \in \binom{[n]}d$.  Let $x^\alpha_X$ be the 
greatest  positive 
integer less than $a_1$ such that $F^\alpha_X(x^\alpha_X)=-1$,
and set $x^\alpha_X=0$ if there is no such integer.  Similarly, set 
$z^\alpha_X$ to be the least integer less than or equal to $n$ such that 
$F^\alpha_X(z^\alpha_X)=-1$, and set $z^\alpha_X=n+1$ if there is no such 
integer.  

\proclaim{Definition 2}
For $d=0$, define $f_2$ as in Definition 1.  

For $d>0$, 
let $\alpha \in B(n,d)$. Let $X = \{a_1,\dots,a_{d}\}_< \in 
\binom{[n]}{d}$.  If, $\forall y \not \in X$ such that $a_1<y<a_d$, 
$F_X^\alpha(y)=1$, then 
we associate to $X$ a  simplex $\{x^\alpha_X,a_1,\dots,a_d,z^\alpha_X\}$.
Define $f_2(\alpha)$ 
to be the set of simplices 
associated to some $X$.  \endproclaim

We remark that
it is by no means obvious that $f_2(\alpha)$ forms a triangulation; this
will follow from the following theorem.  

\proclaim{Theorem 4.1} The map $f_2$ satisfies Definition 1.  More specifically,
fix $\pi \in A(n,d)$.  Let $Y_1,\dots,Y_{\binom nd}$ be 
the elements of $\binom{[n]}d$ under the order $\pi$.  Define $\alpha_i$ and
$S_i$ as in Definition 1.  Then if $S_i =S_{i+1}$ then $f_2$ associates no
simplex to $Y_{i+1}$, and if $S_i$ and $S_{i+1}$ do not coincide, then
$f_2$ does associate a simplex to $Y_{i+1}$, and $S_i$ and $S_{i+1}$ differ
in that $S_i$ contains the bottom facets of this simplex and $S_{i+1}$ 
contains its top facets.  
\endproclaim

\demo{Proof} The proof is by induction on $d$.  
The assertion is clear for $d=0$.  It is also straightforward to check for
$d=1$.  So assume $d>1$.

Fix some $0 \leq i < \binom nd$.  For simplicity, we will denote $Y_i$ by
$Y$.  
Let $Y=\{ a_1,\dots, a_d\}_<$.  For $1\leq k \leq d$, let $L_k$ be the 
$d-1$-dimensional face 
$$L_k(j)=\left\{ \matrix F^{[\pi]}_Y(j) & j \ne a_k\\ (-1)^{k+d+1} & j = a_k 
\endmatrix \right. $$
Let $T_k$ be the face:
  $$T_k(j)=\left\{ \matrix F^{[\pi]}_Y(j) & i \ne a_k\\ (-1)^{k+d} & j = a_k 
\endmatrix \right. $$
By Lemma 2.1, 
the $L_k$ are the lower faces of $F^{[\pi]}_Y$, and $T_k$ are its upper
faces.  Thus, the $L_k$ are in $K(\alpha_i)$, while the $T_k$ are in 
$K(\alpha_{i+1})$.   

We now calculate 
the images of the $L_k$ and $T_k$ under $f_2$.  Suppose first that
there is some $j \not \in Y$, $a_2<j<a_{d-1}$, such that $F_Y^{[\pi]}(j)=-1$.
In this case, $F_Y^{[\pi]}$ has no simplex associated to it, and
none of the $L_k$ or $T_k$ have a simplex associated to them.  
Thus, it follows that $S_i=S_{i+1}$.  So Definition 1 says there should be
no simplex associated to $Y$.  And in this Definition 2 concurs.  

Suppose next that there is some $a_1<j<a_2$ and some $a_{d-1}<j'<a_d$, such
that $F_Y^{[\pi]}(j)=-1$ and $F_Y^{[\pi]}(j')=-1$.  In this case, the same
thing happens: none of the faces have simplices associated to them.  

Suppose now that there is no $a_{d-1}<j'<a_d$ with $F_Y^{[\pi]}(j')=-1$,
but there is some
$a_1<j<a_2$ such that $F_Y^{[\pi]}(j)=-1$.  Let $j$ be the greatest such.  
Then the only one among
the $L_k$ which has a simplex associated to it is $L_1$, and the simplex
associated to it is $\{j,a_2,\dots,a_d,z^{[\pi]}_Y\}$.  
Similarly, the only one 
among
the $T_k$ which has a simplex associated to it is $T_1$, and again, the 
simplex associated to it is $\{j,a_2,\dots,a_d,z^{[\pi]}_Y\}$.  Thus, Definition 
1 
says there should be no simplex associated to $Y$, and Definition 2 concurs.

The case where there is no $a_1<j<a_2$ with $F_Y^{[\pi]}(j)=-1$, but
there is some $a_{d-1}<j'<a_d$ with $F_Y^{[\pi]}(j')=-1$, is dealt with
in exactly the same way.  

Finally, we consider the case where there is no $j \not \in Y$, $a_1<j<a_d$,
such that $F_Y^{[\pi]}(j)=-1$.  
Then there is a simplex associated to $L_k$ for all $k$ of the opposite parity
to $d$, to $k=d$, and to $k=1$ regardless of the parity of $d$.  Similarly,
there is a simplex associated to $T_k$ for all $k$ of the same parity as 
$d$, and also for $k=1$ regardless of the parity of $d$.  
It is straightforward to check that
the simplices associated to the $L_k$ and those associated to the $T_k$ are
form the bottom and top of the 
simplex $\{x^{[\pi]}_Y,a_1,\dots,a_d,z^{[\pi]}_Y\}$, which is the simplex
associated to $Y$.  This completes the
proof of the theorem.  
\enddemo

We shall denote by $f$ the map defined by the two equivalent definitions
above.  

\proclaim{Proposition 4.1} The map $f: B(n,d) \rightarrow S([0,n+1],d+1)$ is 
order-preserving. \endproclaim

\demo{Proof}
We will show that if $\beta \gtrdot \gamma$ in $B(n,d)$ then either
$f(\beta)=f(\gamma)$ or $f(\beta)\gtrdot f(\gamma)$.  

It is shown in [MS] that $\0h_d$ is the minimum element of $B(n,d)$, and
$\1h_d$ is the maximum element.  This implies that there is an unrefinable
chain from $\0h_d$ to $\gamma$, and an unrefinable chain from $\beta$ to 
$\1h_d$.  At each step along these chains, the inversion set increases by
a single element.
The order in which these elements are added 
defines an admissible ordering $\pi \in A(n,d+1)$ 
which has the property
that $I(\gamma)$ and $I(\beta)$ occur as initial subsequences.  

Define $\alpha_i \in B(n,d)$ by $I(\alpha_i)=\In_i(\pi)$, and let 
$S_i=f(\alpha_i)$.  Definition 1 of $f([\pi])$ tells us that 
$S_i$ and $S_{i+1}$ either coincide or $S_i \lessdot S_{i+1}$.
There is some $k$ such that $\alpha_k=\gamma$ and $\alpha_{k+1}=\beta$,
so either $f(\gamma)=f(\beta)$ or $f(\gamma) \lessdot f(\beta)$, as desired.
\enddemo

For completeness we give yet another definition of $f$, also from [Ra1].  

\proclaim{Definition 3} Let $\beta \in B(n,d)$.  Let $I(\beta)=\{
X_1,\dots,X_r\}$ where the $X_i$ are ordered so that every initial subsequence
is also consistent.  Set $T_0 = \0h_d \in S([0,n+1],d)$.  
Define $T_i$ by induction: if $T_{i-1}$ contains
the bottom facets of a simplex with vertices $\{x\} \cup X_i \cup \{z\}$ for
some $x$ less than any element of $X_i$ and $z$ greater than any element of
$X_i$, then let $T_i$ consist of the facets of $T_{i-1}$ with these bottom
facets replaced by the simplex's corresponding top facets.  Otherwise, let 
$T_i=T_{i-1}$.  Then set $f_3(\beta)=T_r$.  
\endproclaim

\proclaim{Theorem 4.2 [Ra1]} The map $f_3$ is well-defined and 
coincides with the map $f$.  
\endproclaim

\demo{Proof}
It is shown in [Ra1] that Definition 2 and Definition 3 are equivalent.
We give a somewhat different proof.    
In the proof, we fix $n$ and $d$, and induct on $r$, the size of the 
inversion set of $\beta \in B(n,d)$.  Define $\gamma\in B(n,d)$ by
$I(\gamma)=\{X_1,\dots,X_{r-1}\}$.  By the induction hypothesis, 
$f_2(\gamma)=f_3(\gamma)$.
Now $K(\gamma)$ and $K(\beta)$ differ in that there is some $d+1$-face $F$ 
of the standard $n$-cube such that $K(\beta)$ includes the top facets and
$K(\gamma)$ the bottom facets.

As in the proof of Proposition 4.1, we construct an order $\pi \in A(n,d+1)$,
$\alpha_i \in B(n,d)$ such that $I(\alpha_i)=\In_i(\pi)$, so that there
exists some $k$ such that $\alpha_k=\gamma$ and $\alpha_{k+1}=\beta$.  

Now, as shown in the proof of Theorem 4.1, exactly one of three things can 
happen:

i) The map $f_2:B(n,d)\rightarrow S([0,n+1],d+1)$ doesn't associate a simplex
to any of the facets of $F$, and $f_2(\gamma)=f_2(\beta)$.

ii) The map $f_2$ as above associates a simplex to exactly one upper facet of
$F$ and one lower facet of $F$ and $f_2(\gamma)=f_2(\beta)$.

iii) We have $f_2(\gamma)\lessdot f_2(\beta)$.  

In cases (i) and (ii), it is easy to check that $f_3(\gamma)$ does not 
contain the bottom facets of a $d+1$-simplex with vertices $\{x\} \cup
X_r \cup \{z\}$ as above, so $f_3(\beta)=f_3(\gamma)=f_2(\gamma)=f_2(\beta)$, 
as desired.

In case (iii), $f_3(\gamma)$ contains the bottom facets of a simplex of the
desired form, and $f_2(\beta)$ consists of $f_2(\gamma)$ with these bottom
facets replaced by the corresponding top facets.  Thus $f_3(\beta)=f_2(\beta)$,
as desired.  \enddemo

It is claimed (without proof) in [KV] that $f:B(n,d) \rightarrow
S([0,n+1],d+1)$ is surjective for all $n$ and $d$.  We cannot prove this
in general.  In the following two sections, we will consider the cases
$d=1$, where surjectivity will turn out to be equivalent to known results,
and $d=2$, where surjectivity is new.  

\head 5. The map $B(n,1) \rightarrow S([0,n+1],2)$ \endhead

In this section we show that the map $B(n,1) \rightarrow S([0,n+1],d)$
is essentially the same as a very familiar map from permutations to planar
binary trees (see [St1, 1.3.13],[BW1],[LR1],[To]).  Because this map 
appears in many guises,
we will give our own definition which is equivalent to all the others.

Let $Y_n$ denote the planar binary trees with $n$ internal vertices.
For $(a_1,\dots,a_p)$ a sequence of $p$ distinct numbers, let $\std(a_1,\dots,a_p)$, the standardization of $(a_1,\dots,a_p)$, denote the sequence of 
numbers from 1 to $p$ arranged in the same order.    
Define $\psi:B(n,1) \rightarrow Y_n$ inductively, as follows: for $n$=0,
$\psi$ applied to the empty permutation is an empty tree; and for 
$n\geq 1$, $\pi \in B(n,1)$, write $\pi = ( a_1 \ \dots \ a_p\  n \ b_1 \ 
\dots \ b_q)$, and then let $\psi(\pi)$ be the tree consisting of one
parent node with left subtree $\psi(\std(a_1,\dots,a_p))$, and right
subtree $\psi(\std(b_1,\dots,b_q))$.  

We recover the map $f:B(n,1) \rightarrow S([0,n+1],2)$ by composing 
$\psi$ with a standard bijection $\theta$ between 
triangulations of an $n+2$-gon
and planar binary trees with $n$ internal vertices, as follows.  
Choose a 
geometric realization of $C([0,n+1],2)$, which we will refer to as $P$.
Fix $S\in S([0,n+1],2)$.  $S$ can be viewed as a triangulation of $P$.
Put a vertex inside each triangle of $S$.  Connect two vertices if their 
triangles share a common edge.  Orient the edge joining the vertices so 
that it points from the
triangle above the edge to the one below the edge (``above'' and ``below''
are with respect to the second coordinate).  For each external edge of $P$, 
except that between 0 and $n+1$, attach a leaf to the vertex corresponding
to the triangle containing that edge.  Because of the way we drew $P$, every
triangle has one upper edge and two lower edges, which we may view as a left
edge and a right edge.  Thus, the tree we have drawn can be viewed as a 
planar binary tree, whose root is the vertex associated to the triangle 
containing the edge from $0$ to $n+1$.  Let this planar binary tree be 
denoted $\theta(S)$.  In the diagram below, we see a triangulation of a 
$5$-gon and its corresponding tree superimposed.  
It is clear that $\theta$ is a bijection.  

$$\epsfbox{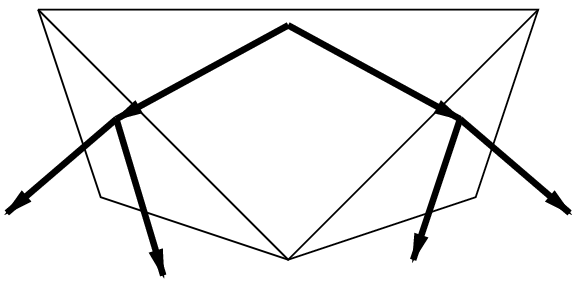}$$
%Example diagram

\proclaim{Proposition 5.1} For $\pi \in B(n,1)$, 
$f(\pi)=\theta^{-1}(\psi(\pi))$. \endproclaim

\demo{Proof} The proof is a simple inductive check.  \enddemo

We now wish to describe the fibre of $f$ over $S \in S([0,n+1],2)$.  These
results have already appeared in the literature; see 
[LR1], [LR2], [BW1], [BW2].  

Any triangle in $S$ has a unique middle vertex, the vertex between
the two bottom edges.  
If the vertices of the triangle are $a<b<c$, the
middle vertex is $b$.  Each vertex of $P$ other than $0$ and
$n+1$ is the middle vertex of a unique triangle of $S$: the triangle
with $b$ as a middle vertex is the one containing $b$ and points immediately
above it.  

We move briefly into greater generality.  Let $T \in S([0,n+1],d)$.  A 
linear order on its simplices is said to be {\it ascending} if for any pair of
simplices sharing a facet, the simplex lying above the intersection facet
follows  the simplex below the intersection.  
It is shown in [Ra1] that there exist ascending  orders on the simplices
of any triangulation of a cyclic polytope of arbitrary dimension.

A linear order on the triangles of $S$ corresponds to a permutation of
$[n]$ by mapping triangles to middle vertices.  

\proclaim{Proposition 5.2} The fibre of $f$ over $S$ consists of the 
permutations 
corresponding to ascending  orders on the triangles of $S$.  \endproclaim

\demo{Proof} The proof is immediate from Definition 1.  \enddemo

This motivates us to inquire further about the ascending  orders on
the triangles of $S$.  We see that the final triangle must be the one
containing the edge $\{0,n+1\}$, which we shall denote $A$. Preceding it
must be a shuffle of an
ascending order on the triangles to the left of $A$ and an ascending order
on the triangles to the right of $A$.  This allows us to prove the following
(already known) proposition:

\proclaim{Proposition 5.3} There are maps $\Min, \Max: S([0,n+1],2)\rightarrow 
B(n,1)$
such that 
the fibre of $f$ over $S \in S([0,n+1],2)$ is the non-empty closed interval
$[\Min(S),\Max(S)]$ in $B(n,1)$.  \endproclaim

\demo{Proof}
Let $A$ be the triangle of $S$ containing the edge $\{0,n+1\}$, and let
its bottom vertex be $a$.  
Let $P_l$ and $P_r$ be the subpolygons of $P$ to the left and right
respectively of $A$.  Let $S_l$ and $S_r$ be the restrictions of 
$S$ to $P_l$ and $P_r$ respectively.  By induction, there is a permutation of
$\Min(S_l)$ of $\{1,\dots,a-1\}$ which is the mininimum among those ascending
 with respect
to $S_l$, and similarly a permutation $\Min(S_r)$ of 
$\{a+1,\dots,n\}$ which is the
minimum among those 
ascending with respect to $S_r$, and similarly permutation $\Max(S_l)$ and
$\Max(S_r)$.  Now it is clear that the minimum ascending order with
respect to $S$ is $ ( \Min(S_l) \ \Min(S_r)\ a)$, and the maximum ascending
order with respect to $S$ is $(\Max(S_r)\  \Max(S_l) \ a)$.   
\enddemo

\head 6. The map $B(n,2) \rightarrow S([0,n+1],3)$ \endhead

This section is chiefly dedicated to the proof of the following proposition:

\proclaim{Proposition 6.1} $f:B(n,2)\rightarrow S([0,n+1],3)$ is surjective.
\endproclaim

\demo{Proof}
Let $S\in S([0,n+1],3)$.  Fix an ascending order (as defined in the previous
section) on its simplices: 
$A_1,\dots,A_s$.  

We will now define a chain of $T_i$ in $S([0,n+1],2)$,
having the property that their simplices occur as facets of the simplices
of $S$.  
Let $T_0$ be the minimum element of $S([0,n+1],2)$.
Define $T_i$ inductively by replacing the
simplices of $T_{i-1}$ which are bottom facets of $A_i$ by the top facets
of $A_i$.  

\proclaim{Lemma 6.1} For $i<j$, $\Min(T_i)<\Min(T_j)$ in $B(n,1)$.  
\endproclaim

\demo{Proof} Clearly, it suffices to consider the case where $j=i+1$.  
Suppose the vertices of $A_{i+1}$ are $a<b<c<d$.  Then $T_i$ and $T_{i+1}$
look like:

\smallskip

$$\epsfbox{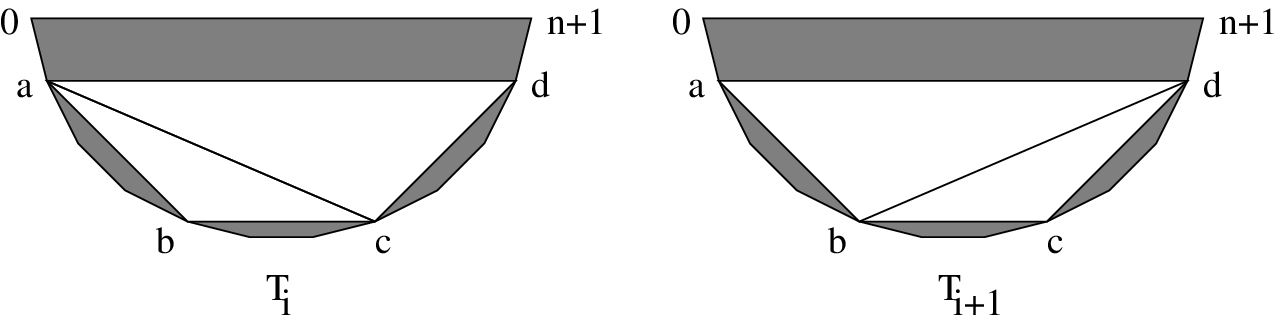}$$
 
\smallskip

Outside the quadrilateral with vertices $a$, $b$, $c$, and $d$, 
$T_i$ and $T_{i+1}$ 
coincide.  

Let $U$, $V$, and $W$ denote their common restrictions to the regions below the
lines $\{a,b\}$, $\{b,c\}$, and $\{c,d\}$ respectively.  Then 
$\Min(T_i)$ and $\Min(T_{i+1})$ coincide except for a consecutive sequence
describing the triangles below the line $\{a,d\}$, which runs
$\Min(U) \ \Min(V) \ b \ \Min(W) \ c $ in $T_i$ and 
$\Min(U) \ \Min(V)\  \Min(W) \ c \ b$ in $T_{i+1}$.  From this it follows
that $\Min(T_{i+1})> \Min(T_{i})$ in $B(n,1)$.  \enddemo

From the lemma, the proposition is almost immediate.  
Pick a maximal chain $\pi_0<\pi_1<\dots<\pi_{\binom n2}$ in 
$B(n,1)$ which refines the chain 
$\Min(T_0)<\Min(T_1)<\dots<\Min(T_s)$.  So $I(\pi_{i+1})$ consists of one more 
element 
than $I(\pi_i)$.  Let $Y_i=I(\pi_{i+1})\setminus I(\pi_i)$.  
Now $Y_1,Y_2,\dots,Y_{\binom{n}{2}}$ is an 
admissible order on $\binom{[n]}{2}$, and thus defines an element of 
$\alpha \in B(n,2)$.  
It is clear from Definition 1 that $f(\alpha)=S$.  

\enddemo

Having shown surjectivity of $f:B(n,2)\rightarrow S([0,n+1],3)$, we might
hope to prove an analogue of Proposition 5.3.  However,
the fibres of the map from $B(n,2)$ to $S([0,n+1],3)$  are more complicated 
than the fibres of the map from $B(n,1)$ to $S([0,n+1],2)$,
as the following example
shows.  Let $I=\{123,124,456,356\}\subset \binom{[6]}2$.  
Let $\alpha\in B(6,2)$ be defined by
$I(\alpha)=I$. Then
$S=f(\alpha)$ is the triangulation consisting of the following simplices:
$$S=
\{0125,0156,0167,0234,0245,1256,1267,2345,2357,2567,3457\}.$$

Using Definition 2, one checks that 
$\beta \in f^{-1}(S)$ iff $I \subset I(\beta) 
\subset I \cup \{134,346\}$.  Now $I\cup \{134\}$ and $I \cup \{346\}$ are
consistent, but $I \cup \{134,346\}$ is not.  Thus, the fibre of $f$ over
$S$ has no maximum element.  
(This example was based on an example
given in [Zi] to show that $B(6,2)$ is not a lattice.)  Thus, there is no
simple analogue of Proposition 5.3.  In turn, this makes it harder to understand
the map $B(n,3)\rightarrow S([0,n+1],4)$.  

\head 7. Interpretation of $\lk_0 \circ f$ and $\lk_{\{0,n+1\}} \circ f$ \endhead

Let $S \in S([0,n+1],d+1)$.  Then $\lk_0(S)=
\{A \setminus \{0\} \mid 0\in A \in S\}$,
the link of $S$ at 0.  
As was remarked in [ER] and proved in [Ra1], thinking of this link as 
describing faces of the vertex figure of $\Delta_{n+1}$ at 0, we see that
$\lk_0(S)\in S(n+1,d)$.  Similarly, we can define $\lk_{n+1}(S)\in 
S([0,n],d)$.
The map $\lk_0$ is order-preserving; $\lk_{n+1}$ is order-reversing.  
(One might
wonder about taking links at other vertices.  For $d$ even, these other
links are not naturally elements of $S(n+1,d)$; for $d$ odd one
can define a link in a suitably labelled $S(n+1,d)$ but this link 
map does not
respect the poset structures.)

In this section, we give geometric interpretations of $\lk_0 \circ f$ and
$\lk_{\{0,n+1\}}\circ f$.  

\proclaim{Proposition 7.1} Let $\pi \in B(n,d)$.  Then a  simplex $A \in 
\lk_{\{0,n+1\}}\circ f$ iff $F^\pi_A$ contains $(1,\dots,1)$.  In other
words, $\lk_{\{0,n+1\}}(f(\pi))$ is the vertex figure of $K(\pi)$ at
$(1,\dots,1)$.\endproclaim

\demo{Proof} The first statement follows immediately from Definition 2 of $f$.
The second statement follows immediately from the first. \enddemo

Let us write $[-1,1]^n_\leq$ for the weakly increasing $n$-tuples from
$[-1,1]$.    This set forms a simplex with $n+1$ vertices, whose coordinates 
consisting of a string of $-1$s followed by a string of $1$s.  We identify 
this simplex with
the standard $n$-simplex by labelling the vertex 
whose first $a-1$ coefficients are $-1$ as vertex $a$.

Let us define a map: 
$$
\align
&W:[-1,1]^n\rightarrow [-1,1]^n_\leq \\
&W(a_1,\dots,a_n)=(a_1,\max(a_1,a_2),\dots,\max(a_1,a_2,\dots,a_n))
\endalign$$

\proclaim{Proposition 7.2} Let $\pi \in B(n,d)$.  Then a simplex 
$A \in \lk_0(f(\pi))$ iff $\dim(W(F^\pi_A))=\dim(F^\pi_A)$.  Consequently,  
$\lk_0(f(\pi))=W(K(\pi))$. 
\endproclaim

\demo{Proof} The first statement follows as before from Definition 2 of $f$.

From the first statement, 
it follows that $\lk_0(f(\pi))$ consists of the images under $W$ 
of 
the faces of $K(\pi)$ which don't drop dimension under $W$.  
Next, we check that if $F^\pi_A$ drops 
dimension, then $W(F^\pi_A)$ coincides with the images under $W$ of the 
union of the upper faces of $F^\pi_A$, and also with the 
union of the lower faces of $F^\pi_A$, which is straightforward.  Now let 
$A_1,\dots,A_\binom{n}{d}$ be an admissible order in the 
equivalence class $\pi$.  Now let 
$J \subset \left[ \binom {n}d\right]$ 
denote the indices $j$ such that $F^\pi_{A_j}$ does not drop dimension.  
Then a simple induction argument shows that
for any $i$, 
$$W\left(\bigcup_{1\leq j \leq i} F^\pi_{A_j}\right)=
\bigcup_{{1\leq j \leq i\atop j\in J}} W(F^\pi_{A_j})$$
For $i=\binom nd$, this is exactly what we want.  
\enddemo

\head 8. Combinatorics of $S(n,d)$ \endhead

For the remainder of the paper, we shall need a combinatorial description of
$S(n,d)$ introduced in [Th].  We begin with some preliminary definitions.

For $\{a_1,\dots,a_{d+1}\}_<$ a subset of $[n]$, let $r(a_1,\dots,a_{d+1})$
denote the subset of $\binom{[n-1]}d$ which consists of those $d$-sets 
consisting of exactly one element from $[a_i,a_{i+1}-1]$ for  
$1\leq i \leq d$.  Subsets
of $\binom{[n-1]}d$ of this form are called {\it snug rectangles}.
We say that a set of snug rectangles forms a {\it snug partition} if 
each $d$-set in $\binom{[n-1]}d$ occurs in exactly one of the snug 
rectangles.  

To $S \in S(n,d)$, we associate the collection of snug rectangles $r(S)$ 
which consists of the rectanges $r(a_1,\dots,a_{d+1})$ for each simplex
$\{a_1,\dots,a_{d+1}\}_<$ in $S$.  Then we have the following theorem:

\proclaim{Theorem 8.1 [Th]} The map $r$ defines a bijection from $S(n,d)$ to
snug partitions of $\binom{[n-1]}d$.  \endproclaim

The description of the covering relations of $S(n,d)$ in terms of snug 
partitions is straightforward.  As we know, $S \gtrdot T$ in 
$S(n,d)$ is equivalent to the existence of 
some $d+1$-simplex $\{a_1,\dots,a_{d+2}\}_<$ 
such that $S$ and $T$ coincide except within this simplex, where $S$ consists 
of its top facets and $T$ consists of its bottom facets.  By Lemma 3.1, 
this is equivalent to the existence of $\{a_1,\dots,a_{d+2}\}_<$ such that
 $r(S)$ and $r(T)$ coincide except that $r(S)$ contains the
snug rectangles $r(a_1,\dots,\hat a_i,\dots,a_{d+2})$ for $i$ odd,
and $r(T)$ contains the snug rectangles $r(a_1,\dots, \hat a_i,\dots,a_{d+2})$
for $i$ even.  

It is sometimes convenient to adopt a different point of view on snug 
partitions, where we partition $\binom{[n-1]}{n-d-1}$ instead of 
$\binom{[n-1]}{d}$.  For $\{a_1,\dots,a_{d+1}\}_<\subset [n]$, let 
$r^c(a_1,\dots,a_{d+1})$ denote the elements of $\binom{[n-1]}{n-d-1}$
which are complements in $[n-1]$ of an element of $r(a_1,\dots,a_{d+1})$.  
We refer to $r^c(a_1,\dots,a_{d+1})$ as a {\it complementary snug rectangle}.

This complementary snug rectangle can be
described explicitly as follows.  Let $\{a_1^c,\dots,a_{n-d-1}^c\}_<=
[n]\setminus\{a_1,\dots,a_{d+1}\}$.  Then 
$$r^c(a_1,\dots,a_{d+1})= (\{a_1^c-1,a^c_1\}\times \dots\times 
\{a_{n-d-1}^c-1,a_{n-d-1}^c\}) \cap \binom{[n-1]}{n-d-1}.$$

A {\it complementary snug partition} is a partition of $\binom{[n-1]}{n-d-1}$ 
into complementary snug rectangles.  
A complementary snug partition records the same information as a snug 
partition, but sometimes it is handier to deal with.

We now describe an important feature of the combinatorics of $S(n,d)$, namely,
the collapse maps, poset maps from $S(n,d)$ to $S(p,d)$ with $p<n$.  

Let $I$ be a subset of $[n-1]$. 
Let 
$m_I:[n]\rightarrow I\cup \{n\}
$ be the map defined by $m_I(a)=\min\{i\in I\cup\{n\}, i \geq a\}$.

Consider the map from $\R^n$ to $\R^{I\cup\{n\}}$ which takes 
$e_i$ to $e_{m_I(i)}$.  
This defines a map from $\Delta_{n-1}$ to 
$\Delta_{|I|}\subset \R^{I\cup\{n\}}$.  We define a 
map $c_I$ on faces of $\Delta_{n-1}$, which takes a face to its image in
$\Delta_{|I|}$, or to $\emptyset$ if its image is lower dimensional.  
Explicitly, if $A=\{a_1,\dots,a_{d+1}\}_<$, then $$c_I(A)=\{m_I(a_1),\dots,
m_I(a_{d+1})\}$$
provided the $m_I(a_i)$ are all distinct, and $c_I(A)=\emptyset$ otherwise.  

Now, for $S \in S(n,d)$, define $c_I(S)$ to be the collection of non-empty
$c_I(A)$ for $A \in S$.  We have the following lemma:

\proclaim{Lemma 8.1} For $I\subset J \subset [n-1]$, $S\in S(n,d)$,
$A=\{a_1,\dots,a_{d+1}\}_<$:

i) $c_I(S) \in S(I\cup \{n\},d)$

ii) $c_I(S)=c_I(c_J(S))$

iii) $r(c_I(A))=r(A) \cap \binom Id$

iv) $r(c_I(S))= \{X \cap \binom{I}d \mid X \in r(S),
X \cap \binom{I}d \ne \emptyset \}$.
\endproclaim

\demo{Proof}  The first part is perhaps easiest to see if we think of 
$S(n,d)$ as triangulations of $C(n,d)$.  If $I=\{i_1,\dots,i_{r-1}\}$,
we can pick a geometric realization of $C(n,d)$ and gradually deform it,
bringing the vertices $1,\dots,i_1$ closer and closer together, and similarly
for $i_1+1,\dots,i_2$ and so on, while preserving the property of being a 
cyclic polytope. The limit of this process is a cyclic polytope with vertices
labelled by $I\cup\{n\}$.  
If we begin with a triangulation $S$ and deform it in this
manner, discarding simplices which degenerate, we obtain $c_I(S)$. 

The other parts are straightforward.  \enddemo

We now draw some consequences of the fact that collapse maps take 
triangulations to triangulations.  

\proclaim{Lemma 8.2} Let $X=\{x_1,\dots,x_{d+1}\}_<\in \binom{[n-1]}{d+1}$,
and let $P$ be the $d$-packet of its $d$-subsets. Then:

(i) Let $R$ be a snug rectangle in $\binom{[n-1]}{d}$.  
The possible
intersection of $R$ with $P$ are: $\{X \setminus \{x_1\}\}$, $\{X \setminus
\{x_{d+1}\}\}$, or $\{X \setminus \{x_i\}, X\setminus \{x_{i+1}\}\}$ 
for $1\leq i \leq d$.

(ii) Let $S \in S(n,d)$.  
Let $W$ be the corresponding snug partition of $\binom{[n-1]}{d}$.  
Then the non-empty
intersections of $P$ with rectangles of $W$ are either
$$\align &\{X \setminus \{x_{d+1}\},X \setminus \{x_{d}\}\}, 
\{X\setminus \{x_{d-1}\}, X \setminus \{x_{d-2}\}\},\dots 
\text{ or}\\
&\{X\setminus \{x_{d+1}\}\},\{X \setminus \{x_d\}, X\setminus \{x_{d-1}\}\},
\dots
\endalign$$
In the former case $X \in I(S)$; in the latter case $X \not \in I(S)$.  
\endproclaim

\demo{Proof} 
The intersection of a snug rectangle with $P$ must be a snug rectangle in
$\binom Xd$ by Lemma 8.1 (iii).  The snug rectangles in $\binom Xd$
are exactly
the possible intersections listed in the statement of the lemma, which proves
(i).  

The restriction of $W$ to $\binom Xd$ must be a snug partition of $\binom 
Xd$, by Lemma 8.1 (iv).  These correspond to triangulations of 
$C(X,d) \cong C(d+2,d)$; there are two of these.  The snug partition 
corresponding to $\1h_d$ is the first list of rectangles, while the snug
partition corresponding to $\0h_d$ is the second list.  This proves (ii). \enddemo

\head 9. The map $g:S(n,d) \rightarrow B(n-1,d)$ \endhead

In this section we define a map $g:S(n,d) \rightarrow B(n-1,d)$, which
is analogous to $f$ in ways which will be made clear later.  

Observe that $S(d+2,d)$ consists of 2 elements, which, as usual, we denote
$\0h_d$ and $\1h_d$.  For $S\in S(n,d)$, let 
$I(S)=\{ X \in \binom{[n-1]}{d+1} \mid c_X(S)=\1h\}$.  We wish to define 
$g(S)$ by setting $I(g(S))=I(S)$.  In order for this to make sense, we
must prove the following proposition:

\proclaim{Proposition 9.1} For $S\in S(n,d)$, $I(S)\subset \binom{[n-1]}{d+1}$
is a consistent set \endproclaim

\demo{Proof}  To check that a set of $d+1$-subsets of $[n-1]$ is consistent,
we must check that its intersection with any $d+1$-packet is either an 
initial or final segment.  So let $J=\{a_1,\dots,a_{d+2}\}_< \subset [n-1]$.  
By part (ii) of Lemma 8.1, 
$I(S) \cap \binom {J}{d+1}= \{X \in \binom{J}{d+1} | c_X(c_J(S))=\1h\}
=I(c_J(S))$.
Since $c_J(S) \in S(J\cup\{n\},d)\cong S(d+3,d)$, it suffices to consider the proposition
in $S(d+3,d)$.  

The triangulations of $C(d+3,d)$ are well-understood.  There are $d+3$ of them,
and the Hasse diagram $S(d+3,d)$ is shown below.  (The labels will be explained shortly.) 
$$\epsfbox{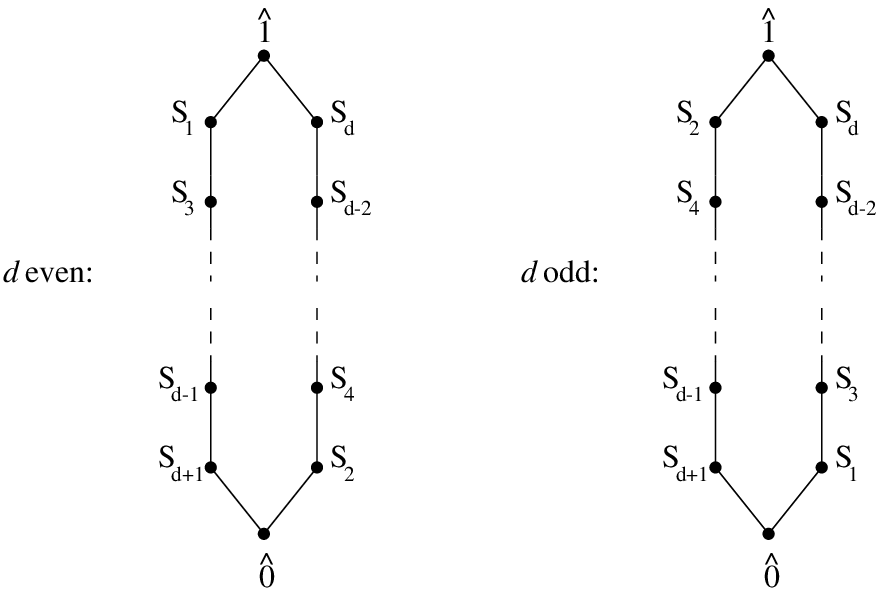}$$

We shall use complementary snug partitions to give a 
 simple pictorial representation of these triangulations.  
Put $\binom{[d+2]}{2}$ in correspondence with a diagram of boxes with
$d+1$ columns numbered 1 to $d+1$ and $d+1$ rows numbered from 2 to $d+2$, 
with boxes in positions $(i,j)$ with $i<j$.  A complementary snug rectangle 
is the intersection of a $2\times 2$ square with 
this diagram, where all partial intersections are allowed, except the 
intersections consisting of a single box $(p,p+1)$ with $1<p<d+1$.  
A complementary snug partition 
consists of a tiling of this diagram by $2\times 2$ squares, with the
stated restrictions on partial intersections.  
It is easy to see by playing with
the pictures that there are two tilings which have no $2\times 2$ square with
its bottom right corner on the diagonal (the left and right pictures below).
The tiling $\0h$ is the one which has the box $(d+1,d+2)$ in a square by 
itself;
the other is $\1h$.  The other tilings have exactly one 
$2 \times 2$ square with its lower right corner on a box on the diagonal
and for 
each box on the diagonal, there is exactly one 
such tiling.  These triangulations are
denoted $S_i$ for $1\leq i \leq d+1$.  $S_i$ 
denotes the triangulation which contains the square 
$\left([i-1,i]\times [i+1,i+2] \right)\cap \binom{d+2}{2}$, whose lower right
corner is $(i,i+1)$.  
$$\epsfbox{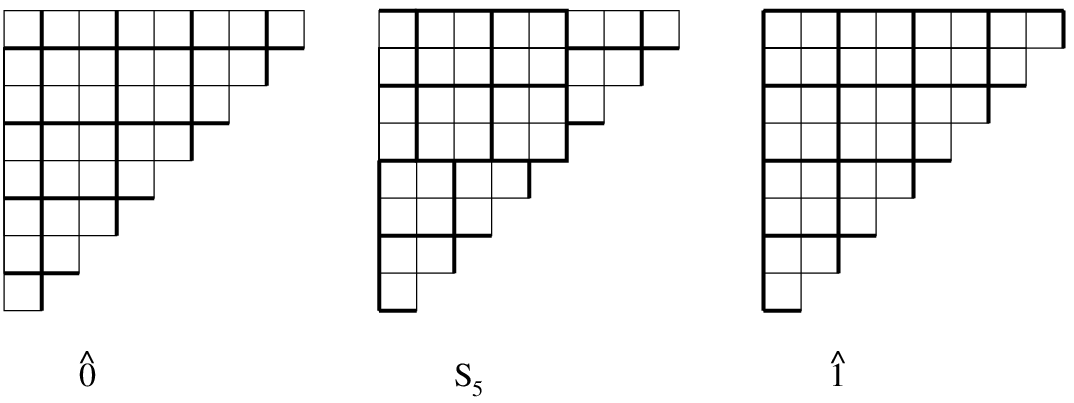}$$

From this analysis, we can see that the Hasse diagram is as above.  Further,
we see that 
$I(\0h)=\emptyset$; $I(\1h)=\binom{[d+2]}{d+1}$; for $d-i$ even
$I(S_i)=\{[d+2]\setminus \{1\},\dots,[d+2]\setminus \{i\}\}$; for $d-i$ odd 
$I(S_i)=\{[d+2]\setminus \{i+1\},\dots,[d+2] \setminus \{d+2\}\}$.
By inspection, all these sets are consistent, which proves the proposition.
\enddemo

In fact, we have proved more than we needed.  
Let us say that $I \subset
\binom{[n-1]}{d+1}$ is {\it superconsistent} if its intersection with any 
$d+1$-packet is either an initial segment of odd length or a final segment 
of the same parity as $d$
(or empty or full).  
Then for $S\in S(n,d)$, $I(S)$ is 
superconsistent.  And more is true:

\proclaim{Theorem 9.1} 
If $I\subset \binom{[n-1]}{d+1}$ is superconsistent, then it is the inversion
set of some $S\in S(n,d)$.  Equivalently, the image of $g:S(n,d)\rightarrow
B(n-1,d)$ consists of those elements with superconsistent inversion sets.
\endproclaim

\demo{Proof} Let $I$ be a superconsistent subset of $\binom{[n-1]}{d+1}$.
We wish to determine a corresponding snug partition of $\binom{[n-1]}{d}$.

For each $X=\{x_1,\dots,x_d\}_< \in \binom{[n-1]}d$, we wish to determine a 
snug rectangle containing $X$. 
Write $\{x_1^c,\dots,x_{n-d-1}^c\}_<$ for $[n-1]\setminus X$.

If $X \cup \{x_i^c\} \in I$, then let 
$$s_i^c=\left\{ \matrix x_i^c &\text{ if $d+i-x_i^c$ is even}\\
                       x_i^c+1 &\text{ if $d+i-x_i^c$ is odd } \endmatrix \right.$$ 

If $X \cup \{x_i^c\} \not \in I$, then let 
$$s_i^c=\left\{ \matrix x_i^c &\text{ if $d+i-x_i^c$ is odd }\\
                       x_i^c+1 &\text{ if $d+i-x_i^c$ is even} 
\endmatrix \right.$$

I claim the $s_i^c$ are all different.  Suppose $s_i^c=s_{i+1}^c$.  Then
clearly $x_i^c=x_{i+1}^c-1$, and either 
$X\cup \{x_i^c\} \in I$, $X \cup \{x_{i+1}^c\} \not \in I$, and $d+i-x_i^c$ is 
odd, or 
$X\cup \{x_i^c\} \not \in I$, $X \cup \{x_{i+1}^c\} \in I$, and $d+i-x_i^c$
is even.  
We assume
that we are in the first case. The number of elements of $X$ greater than
$x_i^c$ is $d+i-x_i^c$, so it is odd.  
Since
$I$ is consistent, if $x\in X$, $x <x_i^c$, then $X\cup\{x_i^c,x_{i+1}^c\}
\setminus \{x\}  \in I$, while if $x \in X$, $x>x_{i+1}^c$, then
$X \cup \{x_i^c,x_{i+1}^c\}
\setminus \{x\} \not \in I$.  But this means that the intersection of $I$ 
with the $d+1$-packet 
$\binom{X \cup \{x_i^c,x_{i+1}^c\}}{d+1}$ is an initial segment of even
length, which contradicts the fact that $I$ is superconsistent.  The 
other case, when $X\cup \{x_i^c\} \not \in I$, is very much the same.  

Now, set $S_X=[n] \setminus \{s^c_1,\dots,s^c_{n-d-1}\}$.  
Now $r^c(S_X)=(\{s_1^c-1,s^c_1\}\times \dots\times 
\{s_{n-d-1}^c-1,s_{n-d-1}^c\}) \cap \binom{[n-1]}{n-d-1}$, and it is clear
that $(x_1^c,\dots,x_{n-d-1}^c)\in r^c(S_X)$, so $X \in r(S_X)$.  

Now, I claim that if $X,Y \in \binom{[n-1]}{d}$, and $Y \in r(S_X)$, then
$S_Y=S_X$.  Let $X=\{x_1,\dots,x_d\}_<$, and let 
$\{x_1^c,\dots,x_{n-d-1}^c\}_<=[n-1]\setminus X$.  
Since we can move from $X^c$ to $Y^c$ by succcessive changes of a single
coordinate by plus or minus one, while remaining in $r^c(X)$, it suffices to
assume that $Y^c=(x^c_1,\dots,x_i^c\pm1,\dots,x_{n-d-1}^c)\in r^c(X)$.  
In fact, we
assume $Y^c=(x^c_1,\dots,x_i^c+1,\dots,x_{n-d-1}^c) \in r^c(X)$.

Thus either $X \cup \{x_i^c \}\in I$ and $d+i-x_i^c$ is odd, or
$X\cup \{x_i^c \} \not \in I$ and $d+i-x_i^c$ is even.  
Let us assume the former; the proof is
similar in the latter case.  

We need to verify
that $X\cup \{x^c_j\} \in I$ iff $Y \cup \{x^c_j\} \in I$ for all
$j \ne i$.  Suppose otherwise for some $j$, and, for the moment, 
suppose $j<i$.  There are an odd number of $x_k$ greater than $x_i^c$,
so, by the superconsistency of $I$, it must be that $X\cup \{x^c_j\}\not \in
I$ while $Y\cup\{x^c_j\}\in I$.   By the consistency of $I$, this
implies that $X\cup\{x^c_i\}\not\in I$, contradicting the assumption.  
The case where $j<i$ is similar.  This completes the proof that
$S_X=S_Y$, and thus that the $r(S_X)$ form a snug partition.  

Let $T$ be the corresponding element of $S(n,d)$.  I claim that
$I(T)=I$.  Let $Z =\{z_1,\dots,z_{d+1}\}_< \subset [n-1]$.  Let
$X=Z \setminus \{z_{d+1}\}$,  $Y=Z \setminus \{z_d\}$.
As usual, let 
$[n-1]\setminus X=\{x_1^c,\dots,x_{n-d-1}^c\}_<$, and similarly for $Y$.  
Let $z_{d+1}=x_i^c$.  Note that $d+i-x^c_i=0$.  Thus, if $Z \not \in I$, 
$s^c_i=x^c_i+1$.  But $y_i^c=z_{d+1}-1$.  Thus $y^c \not \in r^c(S_X)$,
so $X$ and $Y$ are in different snug rectangles in the snug partition
corresponding to $T$.  By Lemma 8.2, this implies that 
$Z \not \in I(T)$.  

Now, suppose that $Z \in I$.  
Observe that for $k\leq z_{d+1}-(z_d+1)$, 
$z_{d+1}-k=x_{i-k}^c$.  As before, $d+i-x^c_i=0$, so $x_i^c=s_i^c$.  
Now, by induction on $k$, using the fact that the $s_j^c$ are all distinct,
it follows that for $k \leq z_{d+1}-(z_d+1)$,
$s_{i-k}^c=x_{i-k}^c$.  
This implies that $y^c \in r^c(S_X)$, so, similarly, $Z\in I(T)$, as 
desired.  This proves the theorem.
\enddemo

Depending on one's taste, one may prefer to justify the construction of 
$S_X$ above by arguing that everything that one has to check to verify
that the $S_X$ form a snug partition 
may be checked after a suitable contraction which reduces the question 
to the case where $n=d+3$, where we have already established the result.  

\head 10. The map $g:S(n,d)\rightarrow B(n-1,d)$ is order-preserving \endhead

\proclaim{Lemma 10.1} Let $S \gtrdot T \in S(n,d)$, and let $S$ and $T$
coincide outside the simplex $\{a_1,\dots,a_{d+2}\}_<$.  Then
$$I(S)=I(T)\adjoin          r(a_1,\dots,a_{d+2})$$
where $\adjoin         $ is disjoint union.   \endproclaim
\demo{Proof} For $I \not \in r(a_1,\dots,a_{d+2})$, one checks that
$c_I(S)=c_I(T)$.

We now investigate $I(S) \cap r(a_1,\dots,a_{d+2})$ and $I(T) 
\cap r(a_1,\dots,a_{d+2})$.   
This is a straightforward check using the definition of $g$,
together with the fact that $S$ contains the snug rectangles 
$r(a_1,\dots,\hat a_i, \dots, a_{d+2})$ for $i$ with the opposite parity to
$d$, while $T$ contains 
those for $i$ having the same parity as $d$.  \enddemo

From this lemma we obtain a corollary which amounts to giving a definition
of $g$ analogous to Definition 3 of $f$, relying as it does to define
$I(S)$ on the choice of a maximal chain from $\0h_d$ to $S$. 

\proclaim{Corollary} Let $S \in S(n,d)$.  Choose an unrefinable chain
$\0h_d=T_0 \lessdot T_1 \lessdot \dots \lessdot T_r=S$.  Let
$R_i$ be the snug rectangle in $\binom {[n-1]}{d+1}$ corresponding 
to the simplex where $T_{i-1}$ and $T_i$ differ.  Then
$$I(S)=\bigcup_{i=1}^r R_i.$$ \endproclaim

We will have occasion to consider a special type of linear order on
the elements of a snug rectangle in $\binom{[n-1]}{d+1}$.  
We say that such an order is 
{\it rectangular} if 

$$\align 
&(x_1,\dots,x_i,\dots,x_{d+1}) > (x_1,\dots,x_i+1,\dots,x_{d+1})
\text { if $(d+1)-i$ is  even} \\
&(x_1,\dots,x_i,\dots,x_{d+1}) < (x_1,\dots,x_i+1,\dots,x_{d+1})
\text { if $(d+1)-i$ is  odd} \tag1\endalign$$

\proclaim{Lemma 10.2} There is a unique rectangular order on
$r(a_1,\dots,a_{d+2})$ up to 
transposition of adjacent pairs of $d+1$-tuples
not both in any $d+1$-packet. \endproclaim

\demo{Proof}
It is  clear that there are rectangular orders.  Suppose that 
$(x_1,\dots,x_{d+1})$ and $(y_1,\dots,y_{d+1})$ are both in 
$r(a_1,\dots,a_{d+2})$ and
lie in some common
$d+1$-packet.  Then, by Lemma 8.2, there is some $i$ such that
$x_j=y_j$ for $j \ne i$.  Thus, it's clear that the order relation
between $(x_1,\dots,x_{d+1})$ and $(y_1,\dots,y_{d+1})$ is determined
by (1) together with transitivity.  
Thus
any two rectangular orders differ only by transpositions of adjacent 
pairs of
elements not both occurring in a common $d+1$-packet.  
\enddemo

We now prove the main result of this section:

\proclaim{Theorem 10.1} The map $g: S(n,d) \rightarrow B(n-1,d)$ is 
order-preserving.  \endproclaim

\demo{Proof}
The statement is clear for $d=0$ and $d=1$, so we may assume that 
$d>1$.  
Let $S\gtrdot T$ in 
$S(n,d)$.  It suffices to show that $g(S)>g(T)$ in $B(n-1,d)$. 

Since $S \gtrdot T$, 
there is some simplex  $A=\{a_1,\dots,a_{d+2}\}_<$ such that
$S$ and $T$ differ in that $S$ contains the top facets of $A$ and $T$ 
contains the bottom facets.  Lemma 10.1 tells us that 
$$I(S)=I(T) \adjoin          r(a_1,\dots,a_{d+2}).$$

The fact that $I(S)$ and $I(T)$ are consistent implies that there
 are admissible orders on $\binom{[n-1]}{d+1}$ which have each of 
$I(T)$  and $I(S)$ as 
initial segments.  We must show that there is a single admissible 
order which has
both as initial segments.  

In fact, we prove something more:

\proclaim{Lemma 10.3} Let $d>1$.  Let $S\gtrdot T \in S(n,d)$.  Let $\alpha$ be 
an order of $I(T)$ such that any initial subsequence is consistent.  
Let $\gamma$ be an order of $\binom{[n-1]}{d+1} \setminus I(S)$ such that
any final subsequence is consistent. Then $\alpha\beta\gamma$ is an admissible
order on $\binom{[n-1]}{d+1}$ iff $\beta$ is rectangular.  
\endproclaim

\demo{Proof} 
We begin by remarking that  there are necessarily
orders $\alpha$ and $\gamma$ as in the statement of the theorem (since 
$\0h_d$ is the unique minimal element and $\1h_d$ the unique maximal element
of $B(n-1,d)$, as shown in [MS]).

Let $\beta$ be a rectangular order on $r(a_1,\dots,a_{d+2})$.
Let us consider the $d+1$-packet $P$ of $d+1$-subsets of $X=\{x_1,\dots,x_{d+2}\}_<$.
We wish to check that it occurs in $\alpha\beta\gamma$ in
either lexicographic order or its opposite.  This is certainly
true if the $d+1$-packet intersects $r(a_1,\dots,a_{d+2})$ in at most
one $d+1$-set.  So suppose it intersects it in more than one place.  
Then by Lemma 8.2, there is some $i$ such that $X \setminus x_i$ and $X\setminus
x_{i+1}$ both lie in $r(a_1,\dots,a_{d+2})$, and this is the entire
intersection of $P$ with the rectangle.  

By the consistency of $I(T)$ and $I(S)$, the intersection
of $I(T)$ and $\binom{[n-1]}{d+1}\setminus I(S)$ with $P$ must be 
$\{X\setminus \{x_j\}\mid j<i\}$ and $\{X\setminus \{x_j\} \mid j>i+1\}$,
but not necessarily respectively.  One now checks that superconsistency of 
$I(T)$ implies that if $d-i$ is odd then $I(T)$ contains the former,
and if $d-i$ is even then $I(T)$ contains the latter, which imply that the 
the elements of $P$ not in $r(a_1,\dots, a_{d+2})$ occur in reverse order if 
$d-i$ is odd, and lex order if $d-i$ is even.  
The rectangularity of $\beta$ ensures that the elements of 
$r(a_1,\dots,a_{d+2}) \cap P$ also occur in the same order.  

On the other hand, if $\beta$ fails to satisfy any of the conditions (1), 
it is clear that there is a $d+1$-packet which intersects 
$r(a_1,\dots,a_{d+2})$ in two places, and these two elements do not occur
in the order which would agree with the order on the rest of the $d+1$-packet. 
Thus, if $\beta$ is not
rectangular, $\alpha\beta\gamma$ is not admissible.  
\enddemo

Theorem 10.1 now follows from Lemma 10.3. \enddemo

\head 11. The map $g:S(n,d) \rightarrow B(n-1,d)$ is a poset embedding \endhead

In this section, we investigate the map $f\circ g$ and show that it
coincides with the ``extension'' map defined in [Ra1].
We then show that $g$ is a poset embedding.

We now
recall Rambau's definition of extension (with a trivial modification
to suit our conventions).  
Let $S\in S(n,d)$.  Then
by definition $\hat S$, the extension of $S$, is 
$$\align \hat S = \{A \cup \{0\} &\mid A \in S\}  \\ 
&\cup \{(x,x+1,a_2,\dots,a_{d+1})
\mid \{a_1,\dots,a_{d+1}\}_<\in A, a_1\leq x\leq  a_2-2\}\endalign$$

It is a nice application of the theory of snug partitions to check 
that $\hat S \in S([0,n],d+1)$.

There is a simple geometrical idea motivating this definition. 
Let $S\in S(n,d)$, thought of as triangulations of $C(n,d)$.  
$S$ defines a hypersurface $\Gamma_S$ in $C(n,d+1)$.  Add a new point on the 
moment curve which precedes all the vertices of $C(n,d)$, and label it
0.  All the faces of $\Gamma_S$ are visible from 0.  $\hat S$ consists
of all the simplices formed by joining $0$ to simplics of $S$, together
with a canonical way to fill in the remainder of $C([0,n],d+1)$.   

It is clear either from this description, or directly from the definition, 
that, as is shown in [Ra1], 
$\lk_0(\hat S)=S$.

\proclaim{Proposition 11.1} For $S \in S(n,d)$, $f (g (S))=\hat S$.  
\endproclaim

\demo{Proof} One checks, using Definition 2 of $f$, that each of the 
simplices of $\hat S$ appears in $f(g (S))$.\enddemo

Since both $f$ and $g$ are order-preserving, 
we recover the result from [Ra1] that the map 
$S \rightarrow \hat S$ is order-preserving.

\proclaim{Theorem 11.1}
The map $g: S(n,d) \rightarrow B(n-1,d)$ is a poset embedding \endproclaim

\demo{Proof} 
Suppose that $S,T \in S(n,d)$, and $g(S)>g(T)$ in $B(n-1,d)$.  Then,
since $f$ is order preserving, $f(g(S))>f(g(T))$.  So
$\hat S > \hat T$, so $S=\lk_0 \hat S >\lk_0 \hat T =T$.  

Thus $g$ is a poset embedding.  \enddemo

\head 12. Alternative definitions of $g$ \endhead

In this section we show that $g$ satisfies two alterative definitions, 
including an analogue of Definition 1 of $f$.  First, we give another
combinatorial construction, for which we need a lemma.  

\proclaim{Lemma 12.1} The unique ascending order on the simplices of the 
triangulation 
$\0h$ of $C(d+2,d)$ is $[d+2]\setminus \{d+2\},[d+2]\setminus \{d\},\dots$.  
The
unique ascending order on the simplices of the 
triangulation $\1h$ is $\dots,[d+2]\setminus \{d-1\},
[d+2]\setminus \{d+1\}$.  \endproclaim

\demo{Proof} This can be seen directly, by examining the intersections of
pairs of simplices in the two triangulations, or by observing that 
an ascending order on 
simplices of a triangulation of $C(d+2,d)$ corresponds to an ascending
chain in $C(d+2,d-1)$, which we have studied in the
proof of Proposition 9.1.  \enddemo

\proclaim{Proposition 12.1} Let $S \in S(n,d)$.  Fix an ascending order on
the simplices of $S$, say, $A_1,\dots,A_r$.  Consider the order on 
$\binom{[n-1]}d$ which consists of the element of $r(A_1)$ followed by
the elements of $r(A_2)$, etc., where the elements within any $r(A_i)$
are written in a rectangular order.  This order is admissible, and
the element of $B(n,d)$ which it defines is $g(S)$.  \endproclaim

\demo{Proof} Let $\pi$ denote an order on $\binom{[n-1]}d$ as in the
statement of the proposition.  
Let $X=\{x_1,\dots,x_{d+1}\}_< \in \binom{[n-1]}{d+1}$.  
Let $P$ denote the $d$-packet of $d$-subsets of $X$.  
Lemma 8.2 describes the two possible sets of non-empty 
intersections of snug rectangles in $r(S)$ with $P$, depending on
whether or not $X \in I(S)$.  
An ascending
order on simplices of $S$ restricts to an ascending order on the simplices 
of $S$ which survive in $c_X(S)$.  Lemma 12.1 describes the unique ascending
order on the simplices of $c_X(S)$.  Thus, the order on $P$ is determined
by Lemma 12.1 and rectangularity, and one checks that this implies that
the elements of $P$ occur in lexicographic order if $X \not \in I(S)$ and
in the reverse of lexicographic order if $X \in I(S)$.  

\enddemo

We now prove the equivalence of a definition of $g$ analogous to the
Definition 1 of $f$.  

\proclaim{Proposition 12.2} Let $d \geq 2$, and $S \in S(n,d)$. Fix an 
ascending order on the simplices of $S$.  Let 
$\0h=T_0 \lessdot T_1 \lessdot \dots \lessdot T_r=\1h$ be the corresponding
chain in $S(n,d-1)$.  Refine the chain $g(\0h) < g(T_1) <\dots<g(\1h)$
to a maximal chain in $B(n-1,d-1)$.  Then $g(S)$ is the element of $B(n-1,d)$
corresponding to that chain.  
\endproclaim

\demo{Proof} Refining the chain as in the statement of the proposition
amounts to finding an admissible ordering on $\binom {[n-1]}{d}$ such that
$I(T_i)$ is an initial subsequence for all $i$.  By Lemma 10.3, there is
a unique element of $B(n-1,d)$ which corresponds to any such order.  This
is the element associated to the admissible order on $\binom{[n-1]}d$
defined in Proposition 12.1, and therefore by that proposition, it 
coincides with $g(S)$.  \enddemo

 Interestingly, this definition fails for $d=0,1$.  Here, different
refinements of the chain of $g(T_i)$ yield different elements of $B(n-1,d)$
(though it is of course easy to specify which refinement to use).

\head 13. Further Directions \endhead

We would like to understand the fibres of $f$ better.  Perhaps, as a 
first step, one might study the fibres of $\lk_0 \circ f$,
since the fibre of $\lk_0 \circ f$ over $S$ has a distinguished element,
namely $g(S)$.  (Contrary to what one might hope, $g(S)$ is neither always
minimal nor always maximal in the fibre.)  

We would also like to 
see the question of the surjectivity of $f$ settled.  

The map $g \circ f:B(n,d) \rightarrow B([0,n],d+1)$ is a map which does 
not seem to have been studied before, and may prove of interest.  

The motivation for [KV] was from the still-developing theory of
$n$-categories.  We hope that our results
may have some application in this area.  In particular, according to
some definitions (see [KV], [St2]), there is an 
$n$-category $\Delta_n$ associated to the $n$-simplex, and an $n$-category
$I_n$ associated to the $n$-cube.  It appears that the map $g$ defines
a map of $n$-categories from $\Delta_n$ to $I_n$ (as the map $f$ was
shown in [KV] to define a map from $I_n$ to $\Delta_{n+1}$).

The order complex of $B(n,d)$ is homotopic to a sphere of dimension
$n-d-2$ [Ra2].  The order complex of $S(n,d)$ is homotopic to a 
sphere of dimension $n-d-3$ [ERR].  Thus, the maps $g$ and $f\circ g$ induce
maps between order complexes which are
homotopy equivalent.  It seems likely that these maps 
are homotopy equivalences.    
(The map $f$ does not induce a map on order complexes because it takes
non-minimal elements to $\0h$.)  

We would also like to understand
the homotopy type of intervals in these posets, or, more restrictedly,
the M\"obius functions of these posets.  There is an 
interesting conjectural description for both, see [Re].  Perhaps the existence
of the new map $g$ will help, at the very least, to connect the
questions for the higher Stasheff-Tamari posets and the higher Bruhat
orders more closely together.  

\head Acknowledgements \endhead

I would like to thank Vic Reiner for comments on an earlier draft
of this paper, and 
Alex Postnikov for suggesting that I investigate the vertex 
figure of an element of 
$B(n,d)$ at $(1,\dots,1) \in [-1,1]^n$.

\enddocument

One topic for further investigation is the matter of the surjectivity of
the map from $B(n,d)$ to $S([0,n+1],d+1)$ for $d>2$, and also the problem
of describing (insofar as possible) the fibres of $f$ for $d>1$.

A recent focus in the study of higher Bruhat orders and higher 
Stasheff-Tamari posets has been the question of the homotopy type of their
order complexes (resolved in [Ra2, ERR]), and more generally of the homotopy
type of intervals in these orders.  (For the relevant conjectures, see [Re].)

Since $f$ generally takes elements 
other than the minimum element of $B(n,d)$ to the minimum element of
$S([0,n+1],d)$, and similarly for maximum elements, $f$ does not induce a 
map on order complexes.  Nonetheless, there are topological questions which
can be asked about $f$ and its fibres, which may shed some light on these
matters.

\enddocument